\documentclass[a4paper,12pt,twoside]{article}
\usepackage[francais]{babel}
\usepackage{amssymb, amsmath, eucal}
\usepackage[all]{xy}
\usepackage{mathrsfs}
\begin{document}
\textwidth15.5cm
\textheight22.5cm
\voffset=-13mm
\newtheorem{The}{Theorem}[subsection]
\newtheorem{Lem}[The]{Lemma}
\newtheorem{Prop}[The]{Proposition}
\newtheorem{Cor}[The]{Corollary}
\newtheorem{Rem}[The]{Remark}
\newtheorem{Titre}[The]{\!\!\!\! }
\newtheorem{Conj}[The]{Conjecture}
\newtheorem{Question}[The]{Question}
\newtheorem{Prob}[The]{Problem}
\newtheorem{Def}[The]{Definition}
\newtheorem{Not}[The]{Notation}
\newtheorem{Ex}[The]{Example}
\newcommand{\C}{\mathbb{C}}
\newcommand{\R}{\mathbb{R}}

\begin{center}

{\Large\bf $L^2$ Extension for Jets of Holomorphic Sections of a Hermitian Line Bundle}

\end{center}

\begin{center}

 {\large Dan Popovici}

\end{center}

\vspace{3ex}

\noindent {\small {\bf Abstract.} Let $(X, \omega)$ be a weakly pseudoconvex K\"ahler manifold, $Y \subset X$ a closed submanifold defined by some holomorphic section of a vector bundle over $X,$ and $L$ a Hermitian line bundle satisfying certain positivity conditions. We prove that for any integer $k\geq 0,$ any section of the jet sheaf $L\otimes {\cal O}_X/{\cal I}_Y^{k+1},$ which satisfies a certain $L^2$ condition, can be extended into a global holomorphic section of $L$ over $X$ whose $L^2$ growth on an arbitrary compact subset of $X$ is under control. In particular, if $Y$ is merely a point, this gives the existence of a global holomorphic function with an $L^2$ norm under control and with prescribed values for all its derivatives up to order $k$ at a point. This result generalizes the $L^2$ extension theorems of Ohsawa-Takegoshi and of Manivel to the case of jets of sections of a line bundle. A technical difficulty is to achieve uniformity in the constant appearing in the final estimate. In this respect, we make use of the exponential map and of a Rauch-type comparison theorem for complete Riemannian manifolds.}

\vspace{3ex}

\subsection{Introduction}\label{subsection:introduction1}

 The problem of extending holomorphic functions, along with an $L^2$ control of the extension, from a subvariety of a complex manifold to this same complex manifold,  was originally solved by T.Ohsawa and K.Takegoshi ([OT87]) and was subsequently generalized by L.Manivel ([Man93]) into the more geometric framework of holomorphic sections of a Hermitian line bundle. The goal of this work is to further generalize the Ohsawa-Takegoshi-Manivel $L^2$ extension theorem to the case of jets of sections of a Hermitian line bundle over a K\"ahler weakly pseudoconvex complex manifold. Specifically, given a complex analytic manifold $X$ of complex dimension $n$, a submanifold $Y$, a Hermitian line bundle $L$ over $X$, and a holomorphic section $f$ of $L$ over $Y$ with good $L^2$ properties, we prove the existence of a holomorphic extension of $f$ to $X$ which satisfies an $L^2$ estimate and has, moreover, locally on $Y$, prescribed partial derivatives up to an arbitrary pre-given order.

Let ${\cal I}_Y$ be the sheaf of germs of holomorphic functions on $X$ which vanish on $Y$. For any integer $k\geq 0$, let ${\cal O}_X/{\cal I}_Y^{k+1}$ be the nonlocally free sheaf of $k$-jets which are ``transversal'' to $Y$. Its fibre at an arbitrary point $y\in Y$ consists of all Taylor series at $y$ truncated to order $k$ along the vertical directions to $Y$. We aim here at extending transversal $k$-jets of sections of the line bundle of holomorphic $L$-valued $(n, \, 0)$-forms, namely sections $f\in H^0(X, \Lambda^n T^{\star}X \otimes L\otimes {\cal O}_X/{\cal I}_Y^{k+1})$. Equivalently, this amounts to extending sections from the unreduced scheme $Y^{(k+1)}$ defined by the quotient sheaf ${\cal O}_X/{\cal I}_Y^{k+1}$ to the ambient manifold $X$.

Assume from now on that the submanifold $Y\subset X$ is defined as \\

\hspace{3ex} $Y= \{ x\in X \, ; \,  s(x) = 0, \, \Lambda^r(ds)(x) \neq 0 \},$ \\

\noindent for some section $s\in H^0(X, \, E),$ generically transverse to the zero section, of a Hermitian holomorphic vector bundle $E$ of rank $r.$ Assume, moreover, $X$ to be equipped with a  K\"ahler metric $\omega$.

The first obstacle to overcome before even stating the result is to define a relevant intrinsic Sobolev-type $L^2_{(k)}$ norm of a $k$-jet. Since the jet sheaf ${\cal O}_X/{\cal I}_Y^{k+1}$ is not locally free, we make the following ad hoc inductive definition. Let $f \in H^0(X, \, \Lambda^n T^{\star}X \otimes L\otimes {\cal O}_X/{\cal I}_Y^{k+1}).$ The holomorphic vector bundle $L':= \Lambda^nT_X^{\star} \otimes L$ is canonically equipped with a metric induced by the metric of $L$ and the reference metric $\omega$ of $X$. Let $\nabla$ be the Chern connection associated with this metric of $L'$, and $\nabla= \nabla^{1, 0} + \nabla^{0, 1}$ its decomposition into its $(1, 0)$ and $(0, 1)$ parts. Fix an arbitrary point $y\in Y$, and let $U$ be a Stein neighbourhood in $X$ giving rise to a surjective morphism $H^0(U, L') \rightarrow H^0(U, L' \otimes {\cal O}_X/{\cal I}_Y^{k+1})$ of local section spaces. Let $\tilde{f} \in H^0(U, L')$ be an arbitrary local lifting of $f.$ Consider now the $C^{\infty}$ vector bundle morphism $T^{\star}X_{|Y} \rightarrow N^{\star}_{Y/X}$ representing the $\omega$-orthogonal $C^{\infty}$ splitting of the exact sequence \\

$0 \rightarrow  N^{\star}_{Y/X} \rightarrow  T^{\star}X_{|Y} \rightarrow T^{\star}Y \rightarrow 0$. \\

\noindent Let $\nabla\tilde{f} = \nabla^{1, 0}\tilde{f} \in H^0(U, L' \otimes T^{\star}X)$. Set $\nabla^1\tilde{f} \in C^{\infty}(U, L' \otimes N^{\star}_{Y/X})$, obtained as a projection of $\nabla^{1, 0}\tilde{f}$ via the surjective bundle morphism $T^{\star}X_{|Y} \rightarrow N^{\star}_{Y/X}$.

\noindent Assume that $\nabla^{j-1}\tilde{f} \in C^{\infty}(U, L' \otimes S^{j-1}N^{\star}_{Y/X})$ has been constructed. Then $\nabla^{1, \, 0}(\nabla^{j-1}\tilde{f}) \in C^{\infty}(U, L' \otimes S^{j-1}N^{\star}_{Y/X} \otimes T^{\star}X)$. We use here the same symbol $\nabla^{1, \, 0}$ to designate the $(1, \, 0)$-type component of the Chern connection on $L' \otimes S^{j-1}N^{\star}_{Y/X}$ equipped with the induced metric. Set $\nabla^j\tilde{f} \in C^{\infty}(U, L' \otimes S^jN^{\star}_{Y/X}),$ the projection of $\nabla^{1, \, 0}(\nabla^{j-1}\tilde{f})$ via the surjective bundle morphisms \\

$S^{j-1}N^{\star}_{Y/X} \otimes T^{\star}X \rightarrow S^{j-1}N^{\star}_{Y/X} \otimes N^{\star}_{Y/X} \rightarrow S^jN^{\star}_{Y/X}$. \\

\noindent  We have thus inductively constructed $\nabla^j\tilde{f} \in C^{\infty}(U, L' \otimes S^j N^{\star}_{Y/X})$ for all nonnegative integers $j$. The pointwise norms $|\tilde{f}|^2 (y), \dots , |\nabla^k\tilde{f}|^2 (y)$ are therefore well defined at every point $y\in Y$ with respect to the metrics canonically induced on the respective vector bundles by the metric of $L'$ and the reference metric $\omega$ of $X$.

\begin{Def}\label{Def:normeponctuelle}

  For a transversal $k$-jet $f\in H^0(U, \, \Lambda^nT^{\star}X \otimes L \otimes {\cal O}_X/{\cal I}_Y^{k+1})$ and a weight function $\rho>0$ on $U$, we define, at every point $y \in Y\cap U,$ the pointwise $\rho$-weighted norm associated to the section $s $ by: \\

$\displaystyle |f|^2_{s, \rho,(k)}(y) := |\tilde{f}|^2(y) + \frac{|\nabla^1\tilde{f}|^2}{|\Lambda^r(ds)|^{2\frac{1}{r}}\, \rho^{2(r+1)}} (y) + \dots + \frac{|\nabla^k\tilde{f}|^2}{|\Lambda^r(ds)|^{2\frac{k}{r}}\, \rho^{2(r+k)}}(y),$ \\

\noindent and the     $L^2_{(k)}$ weighted norm by    : \\

$\displaystyle ||f||^2_{s, \, \rho, \, (k)} = \int\limits_Y |f|^2_{s, \, \rho, \, (k)}\, |\Lambda^r(ds)|^{-2} \, dV_{Y, \, \omega}$.

\end{Def}

\begin{Ex}\label{Ex:exemple}

  {\rm Consider the case where $X=\Omega$ is a bounded pseudoconvex open subset of $\C^n$ containing $0$, $z=(z_1, \dots ,z_n)$ is the coordinate on $\C^n$ and $Y=\{z_1= \dots = z_r =0\}\cap \Omega$. Take $E = \Omega \times \C^r$ equipped with the trivial flat metric, $L = \Omega \times \C$,  and $s = \bigg(\frac{z_1}{e\, \mbox{diam}\, \Omega}, \dots ,\frac{z_r}{e\, \mbox{diam}\, \Omega}\bigg)$. For all $z\in \Omega, |s(z)|^2 = \frac{1}{e^2} \frac{|z_1|^2 + \dots + |z_r|^2}{(\mbox{diam}\, \Omega)^2} \leq \frac{1}{e^2}$. The jet $f$ is then given by holomorphic functions $a_{\alpha}$, $|\alpha| \leq k$, on $Y$, and its weighted $L^2_{(k)}$ norm is given by:

$$\int\limits_Y |f|^2_{s, \, \rho, \, (k)} \, |\Lambda^r(ds)|^{-2} \, dV_{Y, \, \omega} = \int\limits_Y\frac{|a_0|^2}{|\Lambda^r(ds)|^2} \, dV_{Y, \, \omega} + \sum\limits_{|\alpha|=1}\int\limits_Y\frac{|a_{\alpha}|^2}{|\Lambda^r(ds)|^{2\frac{r+1}{r}} \rho^{2(r+1)}} \, dV_{Y, \, \omega} + \dots + $$

$$+ \sum\limits_{|\alpha|=k}\int\limits_Y \frac{1}{(\alpha !)^2} \, \frac{|a_{\alpha}|^2}{|\Lambda^r(ds)|^{2\frac{r+k}{r}} \rho^{2(r+k)}} \, dV_{Y, \, \omega}.$$  }

\end{Ex}

\noindent  It should be noticed that the norm $|f|^2_{s, \, \rho, \, (k)}(y)$ of the $k$-jet $f$ at the point $y\in Y$ is independent of the choice of a local lifting $\tilde{f}$. Indeed, if $\hat{f}\in H^0(U, L')$ is another lifting of $f_{|U}\in H^0(U, L' \otimes {\cal O}_X/{\cal I}_Y^{k+1})$, then $\tilde{f}$ and $\hat{f}$ have the same transversal $k$-jet on $U\cap Y$ (equal to $f_{|U}$). This implies that $\nabla^j\tilde{f} = \nabla^j\hat{f}$ at every point in $U\cap Y$, for all integers $j=0, \dots , \, k$.

\begin{Not}\label{Not:nabla} (a) For a transversal $k$-jet $f\in H^0(U, \Lambda^n T^{\star}_X \otimes L \otimes {\cal O}_X/{\cal I}_Y^{k+1}),$ denote $\nabla^j f:= (\nabla^j \tilde{f})_{|U\cap Y}$, for all $j=0, \dots , \, k$ and an arbitrary lifting $\tilde {f} \in H^0(U, \, \Lambda^n T^{\star}_X \otimes L)$ of $f$.\\

\noindent (b) For every integer $k \geq 0$, set \\

$J^k : H^0(X, \, \Lambda^nT^{\star}_X \otimes L) \rightarrow H^0(X, \, \Lambda^nT^{\star}_X \otimes L \otimes {\cal O}_X/{\cal I}_Y^{k+1})$ \\

\noindent the cohomology group morphism induced by the projection ${\cal O}_X \rightarrow {\cal O}_X/{\cal I}_Y^{k+1}$.

\end{Not}

We can now state the jet extension theorem.

\begin{The}\label{The:principal1}{\bf (Main theorem)} Let $X$ be a complex weakly pseudoconvex manifold of complex dimension $n$, equipped with a K\"ahler metric $\omega$, $L$ a Hermitian holomorphic line bundle, $E$ a Hermitian holomorphic vector bundle of rank $r$ over $X$, and $s\in H^0(X, E)$ a section assumed to be generically transverse to the zero section. Set:

\vspace{1ex}

$Y:= \{ x\in X \, ; \, s(x) = 0, \, \Lambda^r(ds)(x) \neq 0\},$ 

\vspace{1ex}

\noindent  a subvariety of $X$ of codimension $r$. Assume also that, for an integer $k\geq 0,$ the $(1, 1)$-form $i\Theta (L) + (r+k)\, id'd''\log |s|^2$ is semipositive and that there exists a continuous function $ \alpha \geq 1$ on $X$ such that the following two inequalities are satisfied on $X$:

\vspace{1ex}

\noindent (a)\,  $i\Theta (L) + (r+k)\, id'd''\log |s|^2 \geq \alpha^{-1} \, \displaystyle\frac{\{i\Theta (E)s, s\}}{|s|^2}$, \\

\noindent (b)\, $|s| \leq e^{- \alpha}.$ 

\vspace{2ex} 

\noindent If $\Omega \subset X$ is a relatively compact open subset, define an associated weight function $\rho = \rho_{\Omega} >0$ by $\displaystyle \rho(y) = \frac{1}{||D s_y^{-1}|| \, \sup\limits_{\xi\in \Omega}(||D^2s_{\xi}|| + ||Ds_{\xi}||)},$ where $D$ stands for the Chern connection of $E.$ \\

  Then, for every relatively open subset $\Omega \subset X$, and every $k$-jet  $f\in H^0 (X, \, \Lambda^n T^{\star}_X \otimes L \otimes {\cal O}_X/{\cal I}_Y^{k+1})$, such that \\

\hspace{6ex}   $\displaystyle\int_Y |f|^2_{s,\, \rho, \, (k)} \, |\Lambda^r(ds)|^{-2} \, dV_{Y, \, \omega} < +\infty, $ \\

\noindent there exists $F_k \in H^0 (X, \, \Lambda^n T^{\star}_X\otimes L) $ such that $J^k F_k = f $ and \\

$\displaystyle\int_{\Omega} \frac{|F_k|^2}{|s|^{2r} \, (-\log |s|)^2} \, dV_{X, \, \omega} \leq C_r^{(k)} \, \int_Y |f|^2_{s, \, \rho, \, (k)} \, |\Lambda ^r(ds)|^{-2} \, dV_{Y, \, \omega},$ \\

\noindent where $C_r^{(k)} > 0$ is a constant depending only on $r$, $k$, $E$ and  $\sup\limits_{\Omega}||i\Theta(L)||$.  

\end{The}

\noindent {\bf Explanations.} (a) The section $s\in H^0(X, E)$ induces a nowhere zero section $\Lambda^r (ds)$ of the vector bundle $\Lambda^r(T_X/T_Y)^{\star} \otimes \det E$ and its norm $|\Lambda^r(ds)|$ is computed with respect to the induced metric on this vector bundle. The notation $||i\Theta(L)||$ stands for the norm of the curvature tensor of $L$ viewed as a $(1, \, 1)$-form on $X$. It is also worth mentioning that only hypothesis $(a)$ is essential: if $(a)$ holds for a choice of the function $\alpha \geq 1,$ we can always achieve $(b)$ by multiplying the metric of $E$ by a sufficiently small weight $e^{-\chi \circ \psi},$ where $\psi$ is a plurisubharmonic exhaustion of $X$ and $\chi$ is a real convex increasing function. Property $(a)$ still holds after multiplying the metric of $L$ by the weight $e^{-(r+k + \alpha_0^{-1})\chi\circ \psi},$ where $\alpha_0=\inf\limits_{x\in X} \alpha(x).$

(b) Like in the original framework of the Ohsawa-Takegoshi-Manivel extension theorem, it is highly desirable to extend this result to the case of $D''$-closed  differential forms of bidegree $(0, \, q),$ $q\geq 1$. The difficulty stems from the $\bar{\partial}$ operator not being hypoelliptic in bidegree $(0, \, q),$ $q\geq 1.$ We are therefore at a loss for a way of ensuring regularity for the solution. This difficulty, already present in the work of L.Manivel ([Man93]), has not been overcome yet. We refer to ([Dem99], $\S 5$) for details.

(c) The above statement extends straightforwardly to the case where the Hermitian metric of the line bundle $L$ is singular. Indeed, a local singular weight $\varphi$ for such a metric of $L$ can be realized as the decaying limit of a family of $C^{\infty}$ functions $\varphi_{\varepsilon}=\varphi \star \rho_{\varepsilon}$ obtained by convolution with regularizing kernels. Since the constant $C_r^{(k)}$ depends only on the growth rate of the curvature form $i\partial\bar{\partial}\varphi$, a same constant exists for all $\varphi_{\varepsilon}.$ The estimate for $\varphi$, with this same constant, is then obtained by a passage to the limit with $\varepsilon \rightarrow 0.$  \\

The following theorem is a special case of the main theorem for a bounded pseudoconvex open set $\Omega \subset \C^n.$

\begin {The}\label{The:1} Let $\Omega \subset \C^n$ be a bounded pseudoconvex open set and $Y\subset \Omega $ a closed nonsingular subvariety defined by some section $s\in H^0(X, \, E)$ of a Hermitian holomorphic vector bundle $E$ of rank $r$ with bounded curvature form. Assume that $|s| \leq e^{-1}$ on $\Omega$. Let $\rho>0$ be the weight function defined as: \\

$\displaystyle \rho(y) = \frac{1}{||Ds_y^{-1}||\, \sup\limits_{\xi\in \Omega}(||D^2s_{\xi}|| + ||Ds_{\xi}||)},$ \\

\noindent  where $D$ is the Chern connection on $E.$      \\

 Then, for any nonnegative integer $k$ and any plurisubharmonic function $\varphi$ on $\Omega$, there exists a constant $C_r^{(k)} > 0$ depending only on $E$, on $\Omega$ and on the modulus of continuity of $\varphi,$ such that for every holomorphic section $f$ of ${\cal O}_{\Omega}/{\cal I}_Y^{k+1}$ satisfying

\vspace{2ex}

$\displaystyle\int\limits_Y |f|^2_{s, \, \rho, _, (k)}\,  |\Lambda^r(ds)|^{-2} \, e^{-\varphi} \, dV_Y < +\infty$,

\vspace{2ex}

\noindent  there exists a holomorphic function $F_k$ on $\Omega$ such that $J^kF_k = f$ and

$$\int\limits_{\Omega}\frac{|F_k|^2}{|s|^{2r}(-\log |s|)^2} \, e^{-\varphi} \, dV_{\Omega'} \leq C_r^{(k)} \int\limits_Y |f|^2_{s, \, \rho, _, (k)}\,  |\Lambda^r(ds)|^{-2}\, e^{-\varphi} \, dV_Y.$$

\end{The}

The case of a singleton $Y=\{z_0\}$ is particularly interesting. The jet $f$ at $z_0$ is given by complex numbers $a_{\alpha} \in \C$, $|\alpha| \leq k$, $\alpha =(\alpha_1, \dots ,\alpha_n)$. Take $s = (e\, \mbox{diam}\, \Omega)^{-1} \, (z-z_0)$, viewed as a section of the trivial vector bundle $E = \Omega \times \C^n$. It is clear that $|s| \leq e^{-1}$ and that  : \\

\noindent $\displaystyle \int\limits_Y |f|^2_{s, \, \rho, \, (k)}\,  |\Lambda^n(ds)|^{-2} \, e^{-\varphi} =  \Bigg(\sum\limits_{|\alpha| \leq k} |a_{\alpha}|^2 \Bigg) \, e^{-\varphi(z_0)}$. \\

\noindent  Since $-\log |s| = \frac{1}{\varepsilon}\log |s|^{-\varepsilon} \leq \frac{1}{\varepsilon}|s|^{-\varepsilon},$ for all $\varepsilon > 0$, we may replace $|s|^{2n}(-\log |s|)^2$ in the denominator by $|s|^{2(n-\varepsilon)}$. We thus get the following.

\begin{Cor}\label{Cor:2} Let $\Omega \subset \C^n$ be a bounded pseudoconvex open set and let $z_0\in \Omega$ be a point. Then, for every positive integer $k$ and every plurisubharmonic function $\varphi$ on $\Omega,$ there exists a constant  $C_n^{(k)}> 0$ depending only on the modulus of continuity of $\varphi,$ with the following property. For all complex numbers $a_{\alpha}$, $|\alpha| \leq k,$ there exists a holomorphic function $f$ on $\Omega$ such that $f(z_0)= a_0$, $\frac{\partial^{\alpha}f}{\partial z^{\alpha}}(z_0) = a_{\alpha}$, $ 1 \leq |\alpha| \leq k$, and \\

$\displaystyle\int\limits_{\Omega}\frac{|f|^2}{|z-z_0|^{2(n-\varepsilon)}} \, e^{-\varphi(z)} \, dV_{\Omega}(z) \leq \frac{C_n^{(k)}}{\varepsilon^2 \, (\mathrm{diam}\, \Omega)^{2(n-\varepsilon)}} \, \Bigg(\sum\limits_{|\alpha|\leq k}|a_{\alpha}|^2\Bigg) \, e^{-\varphi(z_0)}.$ 

\end{Cor}

We will split the proof of theorem \ref{The:principal1} into two parts. In the first part, the qualitative one, we make use of ideas of the original proof of Ohsawa and Takegoshi ([OT87], [Ohs88]) cast into a more geometric mould by Manivel ([Man93]) and subsequently simplified by Demailly ([Dem99]), that we appropriately fit into our generalized situation. The main idea is to use a ``weight bumping'' technique to concentrate the curvature of the line bundle $L$ on a tubular neighbourhod of the submanifold $Y$. This leads to defining a new curvature operator and to proving $L^2$ estimates for this modified operator which are analogous to those of Hormander. The main tool is a Bochner-Kodaira-Nakano inequality due to Ohsawa ([Ohs95]). This step is performed in section \ref{subsection:demonstrationdutheoremeprincipal} and is common to the proofs of theorems \ref{The:principal1} and \ref{The:1}.  

In the second half of the proof of theorem \ref{The:principal1}, the main goal is to achieve uniformity for the constant appearing in the final $L^2$ estimate. We deal separately with theorems \ref{The:principal1} and \ref{The:1} in proving the quantitative part of the result. In section \ref{subsection:estimationdelasolution}, we apply Cauchy's inequalities to get a control of the growth of the $k$-jet of a holomorphic function in terms of the growth of this very function, and we thus complete the proof of theorem \ref{The:1}. In order to get intrinsic $L^2$ estimates independent of the radius of local holomorphic coordinate patches on $X$ in theorem \ref{The:principal1}, we make use of the exponential map that transfers the situation over to the tangent space to $X$ at a point. In section \ref{subsection:untheoremedecomparaison}, the Jacobi field technique will enable us to get a Riemannian geometric result related to the Rauch comparison theorem. In section \ref{subsection:estimationfinale}, building on this comparison theorem, we get the final estimate in the main theorem thanks to G\aa rding's lemma on the solutions of elliptic systems. Finally, in section \ref{subsection:demonstrationdutheoremevarprincipal}, we dispense with the smoothness restriction on $Y$ through a little standard argument.

\subsection{Ingredients}\label{subsection:rappels}

We list here the main preliminary results underlying the proof of the original Ohsawa-Takegoshi theorem. For proofs and details, see Demailly ([Dem97]).

The main idea in the proof of the Ohsawa-Takegoshi extension theorem ([OT87, Ohs88]) was to derive and use a modified version of the Bochner-Kodaira-Nakano inequality. This version was subsequently improved by Ohsawa ( [Ohs95]) in the following form.

\begin{Prop}{\bf (Main curvature inequality.)}\label{Prop:BKN}

Let $(X, \, \omega)$ be a K\"ahl\'er manifold with a nonnecessarily complete K\"ahler metric, let $(E,h)$ be a Hermitian vector bundle on $X$, and let $\eta, \lambda >0$ be $C^{\infty}$ functions on $X$.

  Then, for every $ u\in {\cal D}(X, \Lambda^{p, \, q}T_X^{\star}\otimes E)$, we have: \\

$||(\eta^{\frac{1}{2}} + \lambda^{\frac{1}{2}})\, D^{''\star}u||^2 + ||\eta^{\frac{1}{2}}\, D''u||^2 + ||\lambda^{\frac{1}{2}}\, D'u||^2 + 2\, ||\lambda^{-\frac{1}{2}} \, d'\eta \wedge u||^2 \geq $\\

\hfill $\geq \langle\!\langle [\eta \, i\Theta(E)-id'd''\eta - i\, \lambda^{-1} \, d'\eta \wedge d''\eta, \Lambda_{\omega}]u, u\rangle\!\rangle$.

\end{Prop}

In the particular case of $(n, q)$-forms, the forms $D'u$ and $d'\eta \wedge u$ vanish as having bidegree $(n+1, q)$. Then the above inequality reads: \\

$||(\eta^{\frac{1}{2}} + \lambda ^{\frac{1}{2}})D^{''\star}u||^2 + ||\eta ^{\frac{1}{2}}D''u||^2 \geq \langle\!\langle [\eta i\Theta (E) -id'd''\eta -i\lambda^{-1}d'\eta \wedge d''\eta, \Lambda]u, u\rangle\!\rangle$.

\vspace{3ex}

 This key curvature inequality enables one to infer the following $L^2$ existence theorem which parallels H\"ormander's $L^2$ existence theorem ([H\"or65, 66]) for a modified curvature operator.

\begin{Prop}\label{Prop:estimations} Let $(X, \, \omega)$ be a K\"ahler manifold. The metric $\omega$ may not be complete but $X$ is assumed to carry a complete K\"ahler metric. Given a Hermitian vector bundle $(E, h)$  and smooth bounded functions $\eta, \lambda >0$ on $X$, consider the curvature operator

\vspace{2ex}

$ B:= B_{E, \omega, \eta, \lambda}^{n, q} := [\eta \, i\Theta(E) - id'd''\eta - i\lambda^{-1} \, d'\eta \wedge d''\eta, \,  \Lambda_{\omega}]$, 

\vspace{2ex}

\noindent acting on the sections of the vector bundle $\Lambda^{n, q}T_X^{\star} \otimes E$, for some $q\geq 1$, and assume that $B$ is positive definite at every point of $X$. \\

  Then, for all $ g\in L^2(X, \, \Lambda^{n, q}T_X^{\star} \otimes E)$ such that $D''g = 0$, and \\

$\displaystyle\int_X\langle B^{-1}g, \, g\rangle \, dV_{\omega} < +\infty,$ \\

\noindent   there exists $f\in L^2(X, \, \Lambda^{n, q-1}T_X^{\star} \otimes E) $ such that $D''f = g$ and \\

$\displaystyle \int_X (\eta + \lambda)^{-1} \, |f|^2  \, dV_{\omega} \leq 2\, \int_X\langle B^{-1}g, \, g\rangle \,  dV_{\omega}$.

\end{Prop}

\vspace{3ex}

 In the course of the proof of the jet extension theorem we shall need to apply the above proposition on a modified metric of the vector bundle under consideration, obtained by multiplying the original smooth metric by the weight $|s|^{-2(r+k)}$ which is singular along $Y=\{s=0\}$. Since the above $L^2$ estimates only work for a smooth metric on a complete K\"ahler manifold, we shall restrict to  $X\backslash Y$. The following standard lemma ensures that $X\backslash Y$ still carries a complete K\"ahler metric.

\begin{Lem}\label{Lem:metkahler}(see, for instance, [Dem82]) Let $(X, \, \omega)$ be a K\"ahler weakly pseudoconvex manifold, $\psi$ a plurisubharmonic exhaustion and $X_c=\{x\in X ; \psi(x) < c \},$ for $c\in \R$. Let $Y=\{ s = 0 \} \subset X$ be an analytic subset defined by a section $s \in H^0(X, E)$ of a Hermitian vector bundle $(E, h)$ over $X$. 

  Then, for all $c\in \R, \hspace{1ex} X_c \backslash Y$ carries a complete K\"ahler metric.

\end{Lem}

\subsection{Proof of theorem \ref{The:principal1}}\label{subsection:demonstrationdutheoremeprincipal}
  
 Assume that the singularity set $\Sigma = \{ s=0, \,  \Lambda^r(ds) = 0 \}$ of $Y$ is empty, which means that $Y$ is a smooth closed subvariety of $X.$ This restriction will be finally lifted by a standard argument in section \ref{subsection:demonstrationdutheoremevarprincipal}. We argue by induction on $k \geq 0$. The case $k=0$ is the d'Ohsawa-Takegoshi theorem. Assume the theorem has been proved for $k-1$. Consider the short exact sequence of sheaves: \\

$ 0 \longrightarrow  S^kN_{Y/X}^{\star} \longrightarrow  {\cal O}_X/{\cal I}_Y^{k+1} \longrightarrow  {\cal O}_X/{\cal I}_Y^k \longrightarrow 0$ \\

\noindent and let $J^{k-1}f \in H^0(X,\, \Lambda^nT^{\star}_X \otimes L \otimes {\cal O}_X/{\cal I}_Y^{k})$ be the image of $f\in H^0(X,\, \Lambda^nT^{\star}_X \otimes L \otimes {\cal O}_X/{\cal I}_Y^{k+1})$ via the induced cohomology group morphism. By induction hypothesis, there exists $F_{k-1} \in H^0 (X, \, \Lambda^nT_X^{\star} \otimes L)$ such that \\

 $ J^{k-1}F_{k-1} = J^{k-1}f $ and 

$$ \int_{\Omega} \frac{|F_{k-1}|^2}{|s|^{2r} (-\log |s|)^2} \, dV_{\omega} \leq C_r^{(k-1)} \, \int_Y |f|^2_{s, \, \rho, \, (k-1)} \, |\Lambda^r (ds)|^{-2} \, dV_{Y, \, \omega},$$ 

\noindent where $C_r^{(k-1)} > 0$ is a constant as in the statement of theorem \ref{The:principal1}. Thus the image of $f-J^kF_{k-1} \in H^0(X,\, \Lambda^nT^{\star}_X \otimes L \otimes {\cal O}_X/{\cal I}_Y^{k+1})$ in $H^0(X,\, \Lambda^nT^{\star}_X \otimes L \otimes {\cal O}_X/{\cal I}_Y^{k})$ is $J^{k-1}f-J^{k-1}F_{k-1} = 0$. This allows the jet $f- J^kF_{k-1}$ to be viewed as a global holomorphic section (on $Y$) of the sheaf $\Lambda ^n T_X^{\star}\otimes L\otimes S^kN_{Y/X}^{\star} = \Lambda ^n T_X^{\star}\otimes L\otimes S^k E^{\star}_{|Y}$.

\vspace{3ex}

\noindent {\bf A $C^{\infty}$ extension of the jet.} We start off by constructing an extension $\tilde{f}\in C^{\infty}(X, \,\Lambda^nT^{\star}_X \otimes L)$ of the holomorphic $k$-jet $f\in H^0(X, \, \Lambda^nT^{\star}_X \otimes L\otimes {\cal O}_X/{\cal I}_Y^{k+1})$ by means of a partition of unity. Consider a covering of $Y$ by coordinate patches $U_i \subset X$ on which the vector bundles $E$ and $\Lambda^nT^{\star}_X \otimes L$ are trivial. Let $e_i$ be a nonvanishing holomorphic section of $\Lambda^nT_X^{\star} \otimes L_{|U_i}$, and $s_1, \dots ,\, s_r$ holomorphic functions on $U_i$ such that $s_{|U_i} = (s_1, \dots , s_r)$ in a trivialization of $E_{|U_i}.$ The functions $s_1, \dots ,\, s_r$ define holomorphic coordinates on $U_i$ transversal to $Y$. Let $z'_{(i)} = (z_{r+1}^{(i)}, \dots ,z_n^{(i)})$ be holomorphic coordinates on $Y\cap U_i$, and write the restriction jet $f$ as $f_{|Y\cap U_i} = w_i\otimes e_{i|Y\cap U_i}$, with $w_i \in H^0(U_i, {\cal O}_X/{\cal I}_Y^{k+1})$. The local $k$-jet $w_i$ is given by holomorphic functions $a_{\alpha}^{(i)}(z'_{(i)})$ on $Y\cap U_i,$ indexed over multi-indices $\alpha=(\alpha_1, \dots , \, \alpha_r) \in \mathbb{N}^r,$ with $|\alpha| \leq k.$ Set \\

$ \hat{f}_i (s, \, z'_{(i)}): = (\sum\limits_{|\alpha|\leq k} a_{\alpha}(z'_{(i)}) \, s^{\alpha} )\otimes e_i \in H^0(U_i, \, \Lambda^nT_X^{\star} \otimes L).$  \\

\noindent Then $\displaystyle \frac{\partial^{\alpha} \hat{f}_i}{\partial s^{\alpha}} (0, \, z'_{(i)}) =  a_{\alpha}^{(i)}(z'_{(i)}),$ for all $\alpha, \, |\alpha| \leq k,$ and $\hat{f}_i$ defines thus a local holomorphic extension of the jet $f$ from $U_i \cap Y$ to $U_i$. Let $\theta_i \in {\cal D}(U_i)$ be a partition of unity such that $\sum \theta_i \equiv 1$ on a neighbourhood of $Y$. Then \\

$\tilde{f}: = \sum\limits_i \theta_i \hat{f}_i  \in C^{\infty} (X, \, \Lambda^nT_X^{\star} \otimes L)$ \\

\noindent defines a $C^{\infty}$ extension of the jet $f$. Furthermore, we have : \\

$\displaystyle D''\tilde{f} = \sum\limits_i d'' \theta_i \wedge \hat{f}_i,$  \hspace{3ex}  $D''\tilde{f} = 0$ \hspace{2ex} on $Y,$  \\

\noindent since all $\hat{f}_i$ assume the same value at every point of $Y$ and $\sum\limits_i d'' \theta_i = 0$ on $Y.$ Likewise, for any multi-index $\alpha =(\alpha_1, \dots , \, \alpha_r) \in \mathbb{N}^r,$ $|\alpha| \leq k,$ if we derive locally $D''\tilde{f}$ along the directions $s=(s_1, \dots , \, s_r)$ transversal to $Y,$ we get: \\

$\displaystyle D^{\alpha}(D''\tilde{f}) = \sum\limits_{\beta \leq \alpha} \sum\limits_i {\alpha \choose \beta} D^{\beta}(d''\theta_i) \wedge D^{\alpha - \beta}\hat{f}_i = 0$ \hspace{2ex} on $Y,$\\

\noindent since for fixed $\alpha - \beta$, all the $D^{\alpha - \beta}\hat{f}_i$ assume the same value at every point of $Y$ (as  $k$-order extensions of the same transversal jet $f$). As the subvariety $Y = \{s=0\}$ is assumed to be smooth, the Taylor development of $D''\tilde{f}$ near $Y$ shows that the $C^{\infty}$extension of $f$ we have just constructed satisfies: \\

  \hspace{6ex}  $|D''\tilde{f}| = O(|s|^{k+1})$   \hspace{3ex} in a neighbourhood of $Y.$

\vspace{2ex}

\noindent {\bf Weight construction; weight bumping technique.} Since we hardly know $\tilde{f}$ away from $Y$, we take a truncation with support in a tubular neighbourhhod of $Y$. Let \\

$\displaystyle G_{\varepsilon}^{(k-1)} := \theta  \bigg(\frac{|s|^2}{\varepsilon ^2}\bigg) \, (\tilde{f} -F_{k-1}) \in C^{\infty}(X, \, \Lambda^nT^{\star}_X \otimes L)$, \\

\noindent where $ \theta : \R \longrightarrow \R$ is a $C^{\infty}$ function such that $\theta \equiv 1$ on $]-\infty, \frac{1}{2}]$, and $ \mbox{Supp}\,\theta \subset { }]-\infty , 1[$. It is clear that $\mbox{Supp}\, G_{\varepsilon}^{(k-1)} \subset \{ |s| < \varepsilon \}$. We shall solve the equation: \\

$(\star) \hspace{2ex} D''u_{\varepsilon} = D''G_{\varepsilon}^{(k-1)}$, \\

\noindent with the extra condition that $\displaystyle \frac{|u_{\varepsilon}|^2}{|s|^{2(r+k)}}\in  L_{loc}^1$ in a neighbourhood of $Y$. This condition ensures that $u_{\varepsilon}$, as well as all its jets of order $\leq k$, vanish on $Y$. Let $\psi$ be a plurisubharmonic exhaustion of $X$, and set $X_c = \{\psi < c\} \subset \subset X$, for all real $c$. The ideal thing would be to solve the equation $(\star)$ on $X$. For technical reasons which will become apparent later, we shall solve the equation $(\star)$ on $X_c \backslash Y_c$, which is still complete K\"ahler thanks to lemma \ref{Lem:metkahler}. The desired holomorphic extension of the jet $f$ will then be $G_{\varepsilon}^{(k-1)} -u_{\varepsilon} + F_{k-1}$. The final solution will be obtained by passing to the limit with $c\rightarrow \infty$ and $\varepsilon \rightarrow 0$.

  Consider now the following functions: \\

 $\displaystyle \sigma_{\varepsilon} : = \log (|s|^2 + \varepsilon ^2),$ $ \eta_{\varepsilon} : = \varepsilon - \chi_0 (\sigma_{\varepsilon})$, $\lambda_{\varepsilon} : = \frac{\chi'_0 (\sigma_{\varepsilon})^2}{\chi''_0 (\sigma_{\varepsilon})}$, \\

\noindent where $\chi_0 : ]-\infty, 0] \rightarrow ]-\infty, 0]$, $\chi_0 (t) = t- \log (1-t)$, for all $t\leq 0$, having the following properties: $\chi (t)\leq t, \hspace{2ex}1 \leq \chi'_0 \leq 2, \hspace{2ex} \chi'' (t) =\frac{1}{(1-t)^2}$. 

\noindent The function $\eta_{\varepsilon}$ is close to $+\infty$ near $Y$ and decays upon getting away from $Y$. It allows therefore to concentrate the curvature of $L$ on a small neighbourhood of $Y$. We define a new curvature operator: \\

$\displaystyle B_{\varepsilon} : = [\eta_{\varepsilon} \, (i\Theta(L) + (r+k) \, id'd'' \log|s|^2) -id'd''\eta_{\varepsilon} - \lambda_{\epsilon}^{-1} \, id'\eta_{\varepsilon} \wedge d''\eta_{\varepsilon}, \, \Lambda]$, \\

\noindent and prove the estimate: \\

$ B_{\varepsilon} \geq \frac{\varepsilon^2}{2|s|^2} (d''\eta_{\varepsilon})(d''\eta_{\varepsilon})^{\star}$, \\

\noindent as operators acting on the $(n, q)$-forms. Easy computations show that \\

 $\displaystyle d'\sigma_{\varepsilon} = \frac{\{D's, \, s\}}{|s|^2 + \varepsilon^2}$, $d''\sigma_{\varepsilon} = \frac{\{s, \, D's\}}{|s|^2 + \varepsilon^2}$, \\  

 $\displaystyle d'd''\sigma_{\varepsilon} = \frac{\{D's, \, D's\}}{|s|^2 + \varepsilon^2} + \frac{\{s, \, D''D's\}}{|s|^2 + \varepsilon^2} - \frac{\{D's, \, s\} \wedge \{s, \, D's\}}{(|s|^2 + \varepsilon^2)^2}$. \\

\noindent   On the other hand, $\Theta(E) = D^2 = D'D'' + D''D'$, and since $D''s = 0,$ owing to $s$ being holomorphic, $D''D's = \Theta(E)s$. This finally yields: \\

$\displaystyle   id'd''\sigma_{\varepsilon} = \frac{i\{D's, D's\}}{|s|^2 + \varepsilon^2} - \frac{i\{D's, s\} \wedge \{s, D's\}}{(|s|^2 + \varepsilon^2)^2} - \frac{\{i\Theta(E)s, s\}}{|s|^2 + \varepsilon^2}$.  \\

\noindent We now use Lagrange's inequality : $\displaystyle  i\{D's, \, D's\} \geq \frac{i\{D's, \, s\} \wedge \{s, \, D's\}}{|s|^2}$ to get : \\

$  id'd''\sigma_{\varepsilon} \geq \frac{\varepsilon^2}{|s|^2} \frac{i\{D's, \, s\} \wedge \{s, \, D's\}}{(|s|^2 + \varepsilon^2)^2} - \frac{\{i\Theta(E)s, \, s\}}{|s|^2 + \varepsilon^2} = \frac{\varepsilon^2}{|s|^2} id'\sigma_{\varepsilon} \wedge d''\sigma_{\varepsilon} - \frac{\{i\Theta(E)s, \, s\}}{|s|^2 + \varepsilon^2}$. \\

\noindent  On the other hand, $ d'\eta_{\varepsilon} = -\chi'_0(\sigma_{\varepsilon}) \, d'\sigma_{\varepsilon}$, $ d''\eta_{\varepsilon} = -\chi'_0(\sigma_{\varepsilon})\, d''\sigma_{\varepsilon}$, and \\

$ 
\begin{array}{lll}
-id'd''\eta_{\varepsilon} &=&  \chi_0'(\sigma_{\varepsilon}) \, id'd''\sigma_{\varepsilon} + \chi''_0 (\sigma_{\varepsilon}) \, id'\sigma_{\varepsilon} \wedge d''\sigma_{\varepsilon} \geq \\

\vspace{3ex}

 & \geq & \displaystyle \bigg(\frac{\varepsilon^2}{2|s|^2} + \frac{\chi''_0(\sigma_{\varepsilon})}{\chi'_0(\sigma_{\varepsilon})^2}\bigg)\, id'\eta_{\varepsilon} \wedge d''\eta_{\varepsilon} -2 \frac{\{i\Theta(E)s, s\}}{|s|^2 + \varepsilon^2}.

\end{array}$

\noindent Let us multiply now the original metric of $L$ by the weight $|s|^{-2(r+k)}$ ; the curvature of this new metric satisfies the inequality \\

$\displaystyle i\Theta(L) + (r+k) \, id'd''\log|s|^2 \geq \alpha^{-1} \, \frac{\{i\Theta(E)s, \, s\}}{|s|^2 + \varepsilon^2},$ \\

\noindent  thanks to hypothesis (a). Indeed, the inequality still holds with the denominator $|s|^2 + \varepsilon^2$ instead of $|s|^2$, owing to the semipositivity of the left-hand term. On the other hand, $|s| \leq e^{-\alpha} \leq e^{-1}$, which entails $\sigma_{\varepsilon} \leq 0$ for $\varepsilon$ small, and \\

$\eta_{\varepsilon} \geq \varepsilon - \sigma_{\varepsilon} \geq \varepsilon - \log(e^{-2\alpha} + \varepsilon^2)$. \\

 \noindent  In addition, we have : $\eta_{\varepsilon} \geq 2\alpha,$ for $\varepsilon < \varepsilon(c)$ small enough. This, along with the previous inequalities, implies: \\

\noindent $ \eta_{\varepsilon}(i\Theta(L) + (r+k)\, id'd''\log|s|^2) - id'd''\eta_{\varepsilon} - \frac{\chi''_0(\sigma_{\varepsilon})}{\chi'_0(\sigma_{\varepsilon})^2}\, id'\eta_{\varepsilon}\wedge d''\eta_{\varepsilon} \geq \frac{\varepsilon^2}{2|s|^2} \, id'\eta_{\varepsilon} \wedge d''\eta_{\varepsilon},$ \\

\noindent on $X_c$. Set $\displaystyle \lambda_{\varepsilon} = \frac{\chi'_0(\sigma_{\varepsilon})^2}{\chi''_0(\sigma_{\varepsilon})}$ and get

\vspace{2ex}

$
\begin{array}{lll}

B_{\varepsilon} : & = & [\eta_{\varepsilon}(i\Theta(L) + (r+k)\, id'd''\log|s|^2) -id'd''\eta_{\varepsilon} - \lambda_{\varepsilon}^{-1} \, id'\eta_{\varepsilon} \wedge d''\eta_{\varepsilon}, \, \Lambda] \geq \\

  & \geq & \displaystyle \bigg [\frac{\varepsilon^2}{2|s|^2}id'\eta_{\varepsilon} \wedge d''\eta_{\varepsilon}, \, \Lambda \bigg ] = \frac{\varepsilon^2}{2|s|^2} (d''\eta_{\varepsilon})(d''\eta_{\varepsilon})^{\star},

\end{array}$

\vspace{2ex}

\noindent as operators acting on the $(n, \, q)$-forms.  \\

\noindent {\bf $\bar{\partial}$-resolution with $L^2$ estimates.} We shall now solve the equation $(\star)$ by means of proposition \ref{Prop:estimations}. In order to avoid the singularity of the weight $|s|^{-2(r+k)}$ along de $Y,$ we prefer working on the relatively compact open subset $X_c \setminus Y_c,$ where $Y_c = Y \cap X_c = Y \cap \{\psi < c\}$, instead of working on $X$ itself. We need verify first that the a priori $L^2$ condition required in proposition \ref{Prop:estimations} is satisfied. Easy computations show that: \\

$D''G_{\varepsilon}^{(k-1)} = g_{\varepsilon}^{(1)} + g_{\varepsilon}^{(2)}$,  \hspace{6ex}  where \\

$ g_{\varepsilon}^{(1)} = (1 + \frac{|s|^2}{\varepsilon^2})\theta'(\frac{|s|^2}{\varepsilon^2}) d''\sigma_{\varepsilon} \wedge (\tilde{f} - F_{k-1}),$ \\

$g_{\varepsilon}^{(2)} = \theta(\frac{|s|^2}{\varepsilon^2}) D''(\tilde{f} - F_{k-1})$.\\

\noindent Since $g_{\varepsilon}^{(2)}$ converges uniformly to $0$ on every compact when $\varepsilon$ tends to $0$, it will have no contribution in the limit. Indeed, $\mathrm{Supp}\, (g_{\varepsilon}^{(2)}) \subset \{|s| < \varepsilon \}$ and $|g_{\varepsilon}^{(2)}| = O(|s|^{k+1}),$ since we have previously shown that $|D''\tilde{f}|  = O(|s|^{k+1})$ in a neighbourhood of $Y.$ This implies that:   \\

\hspace{6ex} $\displaystyle\int\limits_{X_c\backslash Y_c}\langle B_{\varepsilon}^{-1}g_{\varepsilon}^{(2)}, g_{\varepsilon}^{(2)}\rangle \, |s|^{-2(r+k)} \, dV_{X, \, \omega} = O(\varepsilon),$ \\

\noindent if $B_{\varepsilon}$ is locally uniformly bounded below in a neighbourhood of $Y$. If this is not the case, we solve the approximate equation $D''u + \delta ^{\frac{1}{2}} h = g_{\varepsilon}$, where $\delta > 0$ is small (see [Dem99], Remark 3.2, for the details). Since there is no essential extra difficulty in this case, we may assume, for the sake of perspicuity, that we have the desired lower bound for $B_{\varepsilon}$.

\noindent As for $g_{\varepsilon}^{(1)}$, we get the following estimate: \\

$\displaystyle\int\limits_{X_c \backslash Y_c} \langle B_{\varepsilon}^{-1}g_{\varepsilon}^{(1)}, g_{\varepsilon}^{(1)} \rangle \, |s|^{-2(r+k)} \, dV_{X, \, \omega} \leq 8 \int\limits_{X_c \backslash Y_c}|\tilde{f} - F_{k-1}|^2 \theta'\bigg(\frac{|s|^2}{\varepsilon^2}\bigg)^2 \, |s|^{-2(r+k)} \, dV_{X, \, \omega}.$ \\

\noindent  Indeed, $$
                    \begin{array}{lll}

\displaystyle  g_{\varepsilon}^{(1)} & = & -(1 + \frac{|s|^2}{\varepsilon^2})\theta'(\frac{|s|^2}{\varepsilon^2})\chi'_0(\sigma_{\varepsilon})^{-1}\, d''\eta_{\varepsilon} \wedge (\tilde{f} - F_{k-1}), \\

\displaystyle  B_{\varepsilon}^{-1} & \leq & \frac{2|s|^2}{\varepsilon^2}\, (d''\eta_{\varepsilon})^{\star -1}(d''\eta_{\varepsilon})^{-1},

\end{array}$$ and therefore : \\

$\begin{array}{lll}

\displaystyle \langle B_{\varepsilon}^{-1}(d''\eta_{\varepsilon} \wedge u), (d''\eta_{\varepsilon} \wedge u)\rangle &\leq & \frac{2|s|^2}{\varepsilon^2} \langle (d''\eta_{\varepsilon})^{-1 \star} (d''\eta_{\varepsilon})^{-1} (d''\eta_{\varepsilon} \wedge u), (d''\eta_{\varepsilon} \wedge u)\rangle = \\

\displaystyle & = &  \frac{2|s|^2}{\varepsilon^2} \, \langle u, u \rangle = \frac{2|s|^2}{\varepsilon^2} |u|^2.

\end{array}$ \\

\noindent Furthermore, $\frac{2|s|^2}{\varepsilon^2} \leq 2 $ and $(1 + \frac{|s|^2}{\varepsilon^2})\chi'_0(\sigma_{\varepsilon})^{-1} \leq 2,$ on $\mbox{Supp}\,g_{\varepsilon}^{(1)} \subset \{|s| < \varepsilon \}$. This implies \\

$\langle B_{\varepsilon}^{-1} g_{\varepsilon}^{(1)}, \, g_{\varepsilon}^{(1)}\rangle \leq 8 \, \theta'(\frac{|s|^2}{\varepsilon^2})^2 |\tilde{f} - F_{k-1}|^2$. \\

\noindent  If $z=(z_1, \dots , \, z_r)$ is an arbitrary local holomorphic coordinate system transversal to $Y$, we have \\

$\displaystyle \frac{|s|^{2r}}{|\Lambda^r(ds)|^2} = \frac{|z|^{2r}}{|\Lambda^r(dz)|^2}$, \\  

\noindent the norms of the sections $\Lambda^r(ds)\in  H^0(X, \, \Lambda^r(T_X/T_Y)^{\star} \otimes \det E)$ and $ \Lambda^r(dz)\in H^0(U, \, \Lambda^r(T_X/T_Y)^{\star})$ being computed with respect to the metrics induced on the respective vector bundles by $\omega$ and by the given metric on $E$. \\

\noindent The integrand of the last integral can be locally written, after the change of variable $s \leadsto \varepsilon \, s$, as \\

$\displaystyle \frac{|(\tilde{f} - F_{k-1})(\varepsilon \, s, z')|^2}{\varepsilon^{2(r+k)}|s|^{2(r+k)}} \, \frac{\theta'(|s|^2)^2}{|\Lambda^r(ds)|^{2\frac{r+k}{r}}}\, dV_{\omega}(\varepsilon \, s, z') = $\\

$\displaystyle =\frac{|(\tilde{f}-F_{k-1})(\varepsilon \, s, z')|^2}{\varepsilon^{2k}|s|^{2(r+k)}} \, \frac{\theta'(|s|^2)^2}{|\Lambda^r(ds)|^{2\frac{r+k}{r}}}dV_{\omega}(s, z')$. \\

\noindent  Since $J^{k-1}f - J^{k-1}F_{k-1} = 0$, the Taylor series development yields: \\

$\displaystyle (\tilde{f} - F_{k-1})(\varepsilon \, s, z') = \sum\limits_{|\alpha| +|\beta| \geq k} \frac{\varepsilon^{|\alpha| + |\beta|}}{(\alpha + \beta)!} \, \frac{\partial^{\alpha + \beta}(\tilde{f} - F_{k-1})}{\partial s^{\alpha} \partial \bar{s}^{\beta}}(0, z') s^{\alpha} \bar{s}^{\beta} = $ \\

$\displaystyle = \varepsilon^k \bigg(\sum\limits_{|\alpha| = k}\frac{1}{\alpha !}\frac{\partial^{\alpha}(\tilde{f} - F_{k-1})}{\partial s^{\alpha}}(0, z') s^{\alpha} + \sum\limits_{|\alpha| + |\beta| \geq k+1} \frac{\varepsilon^{|\alpha| + |\beta| -k}}{(\alpha + \beta )!} \,  \frac{\partial^{\alpha + \beta}(\tilde{f} - F_{k-1})}{\partial s^{\alpha}\partial \bar{s}^{\beta}}(0, z')  s^{\alpha}\bar{s}^{\beta}\bigg)$= \\

$\displaystyle =  \varepsilon^k (f - J^kF_{k-1})(z') + O(|\varepsilon s|^{k+1}) =  \varepsilon^k \, \nabla^k(f - J^kF_{k-1})(z') + O(|\varepsilon s|^{k+1}). $ \\

\noindent The first sum ranges only on multi-indices $\alpha$ and $\beta$ such that if $|\alpha| + |\beta| = k$, then $|\alpha| = k$. \\

\noindent  This shows that $\displaystyle \frac{|(\tilde{f} - F_{k-1})(\varepsilon s, z')|^2}{\varepsilon ^{2k}}$ converges to $|\nabla^k(f - J^kF_{k-1})(z')|^2$, (see notation \ref{Not:nabla}), uniformly on every compact, when $\varepsilon \rightarrow 0$. \\  

\noindent We have thus proved that: \\

$\displaystyle \int\limits_{X_c \backslash Y_c}\langle B_{\varepsilon}^{-1} g_{\varepsilon}^{(1)}, g_{\varepsilon}^{(1)}\rangle \, |s|^{-2(r+k)} dV_{X, \, \varepsilon} \leq $ \\

$\displaystyle\leq 8 \int\limits_{X_c \backslash Y_c} |\tilde{f} - F_{k-1}|^2 \, \Theta'\bigg(\frac{|s|^2}{\varepsilon^2}\bigg)^2 \, |s|^{-2(r+k)} \, dV_{X, \, \varepsilon} \rightarrow 8 \, C_{r, \, k} \int\limits_{Y_c} \frac{|\nabla^k(f-J^kF_{k-1})|^2}{|\Lambda^r(ds)|^{2\frac{r+k}{r}}} \, dV_{Y, \, \omega},$ \\

\noindent where $\displaystyle C_{r, \, k}:= \int\limits_{z\in \C^r, |z| \leq 1} \theta'(|z|^2)^2 \, \frac{i\Lambda^r(dz)\wedge \Lambda^r(d\bar{z})}{|z|^{2(r+k)}}.$ 

\noindent It is worth noticing that $|\nabla^k(f-J^kF_{k-1})| = |f-J^kF_{k-1}|$, where  $|f - J^k F_{k-1}|$ is the norm of the section: \\

$f - J^kF_{k-1} \in H^0(Y, \, \Lambda^nT^{\star}_X \otimes L \otimes S^kN_{Y/X}^{\star})$ \\

\noindent with respect to the metric induced on $S^kN_{Y/X}^{\star} $ by the reference metric $\omega$ on $X$. Indeed, $S^kN_{Y/X}$ is a subbundle of $(S^kT_X)_{|Y}$ ; we merely take the metric induced on $S^kN_{Y/X}$ by restriction. \\

 The $L^2$ condition required beforehand in proposition \ref{Prop:estimations} is thus satisfied.The solution $u_{c, \varepsilon}$ to the equation $(\star) \,\, D''u_{c, \, \varepsilon} = D''G_{\varepsilon}^{(k+1)} = g_{\varepsilon}^{(1)} +g_{\varepsilon}^{(2)}$ on $X_c \backslash Y_c$ satisfies then the estimate: \\

\noindent $(1)$ \, $\displaystyle \int\limits_{X_c \backslash Y_c}\frac{|u_{c, \varepsilon}|^2}{|s|^{2(r+k)}(-\log (|s|^2 + \varepsilon^2))^2} \, dV_{X, \, \varepsilon} \leq \int\limits_{X_c \backslash Y_c}\frac{|u_{c, \varepsilon}|^2}{(\eta_{\varepsilon} + \lambda_{\varepsilon})|s|^{2(r+k)}} \, dV_{X, \,  \omega} \leq $ \\ 

\noindent $\displaystyle \leq 2 \int\limits_{X_c \backslash Y_c}\langle B_{\varepsilon}^{-1}g_{\varepsilon}, g_{\varepsilon} \rangle |s|^{-2(r+k)} dV_{X, \, \omega} \leq 16C_{r,  k} \int\limits_{Y_c} \frac{|\nabla^k(f - J^kF_{k-1})|^2}{|\Lambda^r(ds)|^{2\frac{r+k}{r}}}  dV_{Y, \, \omega} + O(\varepsilon).$  \\

\noindent  Indeed, we have used the following obvious estimates: \\

$\sigma_{\varepsilon} = \log(|s|^2 + \varepsilon^2) \leq \log(e^{-2\alpha} + \varepsilon^2) \leq -2\alpha + O(\varepsilon^2) \leq -2 + O(\varepsilon^2)$, \\

$ \eta_{\varepsilon} = \varepsilon - \chi_0(\sigma_{\varepsilon}) \leq (1 + O(\varepsilon))\sigma_{\varepsilon}^2$, \\

$\lambda_{\varepsilon} = \frac{\chi'_0(\sigma_{\varepsilon})^2}{\chi''_0(\sigma_{\varepsilon})} = (1 - \sigma_{\varepsilon})^2 + (1 - \sigma_{\varepsilon}) \leq (3 + O(\varepsilon))\sigma_{\varepsilon}^2$, \\

$ \eta_{\varepsilon} + \lambda_{\varepsilon} \leq (4 + O(\varepsilon))\sigma_{\varepsilon}^2 \leq (4 + O(\varepsilon))(-\log(|s|^2 + \varepsilon^2))^2$. \\

\noindent  The extension of $f$ to $X_c \backslash Y_c$ is given by: \\

 $ F_{c, \varepsilon}^{(k)}:= G_{\varepsilon}^{(k-1)} - u_{c, \varepsilon} + F_{k-1}$. \\

\noindent Locally at an arbitrary point of $Y$, this means that all partial derivatives of ordre $\leq k$ of $F_{c, \varepsilon}^{(k)}$ are prescribed by $f$. The function $G_{\varepsilon}^{(k-1)}$ is $C^{\infty}$ on a tubular neighbourhood of $Y$ and $\mbox{Supp}\,G_{\varepsilon}^{(k-1)} \subset \{ |s| < \varepsilon \}$. This implies that: \\

\noindent $(2)$ \, $\displaystyle  \int\limits_{X_c}\frac{|G_{\varepsilon}^{(k-1)}|^2}{(|s|^2 + \varepsilon^2)^r (-\log(|s|^2 + \varepsilon^2))^2} \, dV_{X, \, \omega} \leq \frac{\mbox{Const}}{(\log \varepsilon)^2}.$  \\

\noindent Since \\

$\displaystyle \int\limits_{X_c \backslash Y_c}\frac{|u_{c, \varepsilon}|^2}{|s|^{2r}(-\log(|s|^2 + \varepsilon^2))^2} \, dV_{X, \, \omega} \leq \int\limits_{X_c \backslash Y_c}\frac{|u_{c, \varepsilon}|^2}{|s|^{2(r+k)}(-\log(|s|^2 + \varepsilon^2))^2}\, dV_{X, \, \omega},$ \\

\noindent (1), (2) and the induction hypothesis made on the $L^2$ norm of $F_{k-1}$ imply the estimate: \\

\noindent $\displaystyle \int\limits_{X_c \backslash Y_c}\frac{|F_{c, \varepsilon}^{(k)}|^2}{(|s|^2 + \varepsilon^2)^r (-\log(|s|^2 + \varepsilon^2))^2} \, dV_{X, \, \omega} \leq $ 

$$\leq 16\, C_{r, k} \int_{Y_c}\frac{|\nabla^k(f - J^kF_{k-1})|^2}{|\Lambda^r(ds)|^{2\frac{r+k}{r}}} \, dV_{Y, \, \omega} + \int_{X_c}\frac{|F_{k-1}|^2}{|s|^{2r}(-\log|s|)^2} \, dV_{X, \, \omega} + \frac{\mbox{Const}}{(\log \varepsilon)^2} \leq $$ 

$$\leq 16\, C_{r, k}\int\limits_{Y_c}\frac{|\nabla^k(f - J^kF_{k-1})|^2}{|\Lambda^r(ds)|^{2\frac{r+k}{r}}} \, dV_{Y, \, \omega} + C_r^{(k-1)}  \int\limits_Y|f|^2_{s, \, \rho, \, (k-1)} \, |\Lambda^r(ds)|^{-2} \, dV_{Y, \, \omega} +  \frac{\mbox{Const}}{(\log \varepsilon)^2} $$ 

$$\leq \displaystyle C_r^{'(k)}  \int\limits_Y |f|^2_{s, \, \rho, \, k} \, |\Lambda^r(ds)|^{-2} \, dV_{Y, \, \omega} + 16\, C_{r, k} \int\limits_{Y_c}\frac{|\nabla^k(J^kF_{k-1})|^2}{|\Lambda^r(ds)|^{2\frac{r+k}{r}}} \, dV_{Y, \, \omega} + \frac{\mbox{Const}}{(\log \varepsilon)^2}, $$

\noindent where $C_r^{'(k)} = C_r^{(k-1)} + 16 \, C_{r, \, k}$.

\noindent  We also have $ D''F_{c, \varepsilon}^{(k)} = 0$ on $X_c \backslash Y_c$, by construction. This relation extends from $X_c \backslash Y_c$ to $X_c$ because $F_{c, \varepsilon}^{(k)}$ is $L^2_{loc}$ in a neighbourhood of $Y_c$. This is guaranteed by the following standard lemma on the $\bar{\partial}$ operator (see, for instance, [Dem82]).  

\begin{Lem}\label{Lem:dbarextension}

  Let $\Omega$ be an open subset of $\C^n$ and $Y$ an analytic subset of $\Omega$. Let $v$ be a $(p, q-1)$-form with $L^2_{loc}$ coefficients and $w$ a $(p, q)$-form with $L^1_{loc}$ coefficients such that $d''v = w$ on $\Omega \backslash Y$(in the sense of distributions). Then $d''v = w$ on $\Omega$.

\end{Lem}

\noindent  The ellipticity of the operator $\bar{\partial}$ in bidegree $(0, 0)$ ensures that $u_{c, \varepsilon}$ is $C^{\infty}$. Consequently, $F_{c, \varepsilon}^{(k)}$ is $C^{\infty}$ as well.

We have thus obtained a family of solutions $(F^{(k)}_{c, \varepsilon})_{\varepsilon}$, along with $L^2$ estimates of these, on the relatively compact open subset $X_c$ of $X$. By extracting a weak limit when $\varepsilon \rightarrow 0$, we thus get a solution $F^{(k)}_c$ and an $L^2$ estimate of it on the relatively compact open subset $X_c$, for all $c>0$.

\subsection{Estimation of the solution in theorem \ref{The:1}}\label{subsection:estimationdelasolution}

  In order to get the final estimates in theorems \ref{The:principal1} and \ref{The:1}, it remains to estimate $\displaystyle \int\limits_{Y_c}\frac{|\nabla^k(J^kF_{k-1})|^2}{|\Lambda^r(ds)|^{2\frac{r+k}{r}}} \, dV_{Y, \, \omega}.$

  We will be dealing in this section with theorem \ref{The:1} where the analysis is simplified by the ambient manifold being an open subset $\Omega \subset \C^n.$ We will use the Cauchy inequalities (or, equivalently, Parseval's formula). In the more general case of theorem \ref{The:principal1}, such an approach would yield a constant depending on the radius of the local holomorphic coordinate balls of $X.$ Since this is an uncontrollable quantity, we will avoid this arbitrariness in the subsequent sections by means of the exponential map replacing locally the ambient manifold $X$ by its tangent space at a point. 
  
  Let $\omega$ be the standard K\"ahler metric on $\Omega$. Since the curvature of $E$ is assumed to be bounded, there exists a constant $M > 0$ such that $i\Theta(E) \leq M\omega \otimes \mbox{Id}_E.$ Set $L = \Omega \times \C,$ equipped with the metric of weight $e^{-\varphi - A|z|^2}$, with a constant $A \gg 0$. If we set $\alpha \equiv 1$, the condition $(a)$ of theorem \ref{The:principal1} is equivalent to \\

$\displaystyle id'd''\varphi + A \, id'd''|z|^2 + (r+k)\, id'd''\log |s|^2 \geq \frac{\{i\Theta(E)s, s\}}{|s|^2}$. \\

\noindent Since $id'd''\varphi \geq 0$, $id'd''\log |s|^2 \geq -\frac{\{i\Theta(E)s, s\}}{|s|^2}$, and $ \frac{\{i\Theta(E)s, s\}}{|s|^2}\leq M\omega$, this relation is satisfied as soon as $A$ has been chosen large enough. This choice of $A$ depends on the bound $M$ of the curvature tensor of $E$.

 Let $\psi : \Omega \rightarrow \R$ be a $C^{\infty}$ plurisubharmonic exhaustion of $\Omega$, namely a function such that the level subsets $\Omega_c:= \{ \psi < c \}$ are relatively compact in $\Omega$ for all $c > 0$. We may assume that $\Omega' = \Omega_c$ for some $c,$ and denote $Y_c:= Y \cap \Omega_c$. Consider now a covering of $Y_c$ by open subsets $U_j$, $j=1, \dots , \, p,$ such that on every $U_j$ there exist local coordinates $z=(z', \, z''), \, z'=(z_1, \dots , z_r)$, $z''=(z_{r+1}, \dots ,z_n)$ for which $Y\cap U_j = \{z'=0\}.$ Pick such a $U_j$ and assume that $U_j = B'(0, \, \rho) \times B''(0, \, \rho) \subset B(0, \, \rho \sqrt{2})$, where $B'(0, \, \rho)$ is the ball of radius $\rho$ of $\C^r,$ $B''(0, \, \rho)$ is the ball of radius $\rho$ of $\C^{n-r}$, and $B(0, \, \rho \sqrt{2})$ is the ball of radius $\rho \sqrt{2}$ of $\C^n.$ The jet $\nabla^k(J^kF_{k-1})$ can be written on $U_j$ as $\sum\limits_{|\alpha| = k}\frac{1}{\alpha!} \frac{\partial^{\alpha}F_{k-1}}{\partial z^{'\alpha}}(0, z'') z^{'\alpha}$, and its norm is given by \\

$\displaystyle |\nabla^k(J^kF_{k-1})|^2 = \sum\limits_{|\alpha|=k}\bigg|\frac{\frac{\partial^{\alpha}F_{k-1}}{\partial z^{'\alpha}}(0, z'')}{\alpha !}\bigg|^2 \, e^{-2\varphi (0, z'') - 2A|z''|^2}$.

\noindent  Parseval's formula applied for $z'\in B'(0, \, \rho)$ gives \\

$\begin{array}{lll}

\displaystyle \frac{\mbox{Const}}{\rho^{2r}}\, \int\limits_{z'\in B'(0, \, \rho)} |F_{k-1}(z', z'')|^2 \, d\lambda(z') & = & \displaystyle \sum\limits_{\alpha} \bigg| \frac{\frac{\partial^{\alpha}F_{k-1}}{\partial z^{'\alpha}}(0, z'')}{\alpha !}\bigg|^2 \,  \frac{\rho^{2|\alpha|}}{2r + 2|\alpha|} \geq \\

  & \geq & \displaystyle \sum\limits_{|\alpha|=k}\bigg| \frac{\frac{\partial^{\alpha}F_{k-1}}{\partial z^{'\alpha}}(0, z'')}{\alpha !}\bigg|^2 \, \frac{\rho^{2k}}{2(r+k)}, \\

\end{array}$

\noindent where $\mbox{Const}$ is a universal constant. Consequently, \\

\noindent $ \frac{\displaystyle \sum\limits_{|\alpha|=k}\bigg|\frac{ \frac{ \partial^{\alpha}F_{k-1}}{ \partial z^{'\alpha}}(0, z'')}{ \alpha !}\bigg|^2 \, e^{-2\varphi (0, z'')-2A|z''|^2}}{\displaystyle |\Lambda^r(ds)(0, z'')|^{2\frac{r+k}{r}}} \leq $ \\

\noindent $\displaystyle \leq \mbox{Const} \, \frac{2(r+k)}{\rho^{2(r+k)}}  \int\limits_{z'\in B'(0, \, \rho)} ||F_{k-1}(z', z'')||^2 \, \frac{e^{2(\varphi(z', z'') - \varphi(0, z''))}\, e^{2A|z'|^2}}{|\Lambda^r(ds)(0, z'')|^{2\frac{r+k}{r}}} \, d\lambda(z'),$  \\

\noindent for all $z''\in B''(0, \, \rho),$ where we have denoted by \\

$||F_{k-1}(z', z'')||^2:= |F_{k-1}(z', z'')|^2 \, e^{-2\varphi(z', z'')}\, e^{-2A(|z'|^2 + |z''|^2)},$ \\

\noindent the norm of the section $F_{k-1}$ in the line bundle $L$. Due to an inconsistency in notation, this vector bundle norm $|| \, \, ||$ is the same as the one we had denoted by $|\, \, |$ in the induction hypothesis (see start of section $0. 3$). Let $\varepsilon$ be a modulus of continuity for $\varphi$, namely a function such that \\

$|\varphi(z', z'')-\varphi(0, z'')| \leq \varepsilon(\, |z'|\,), \hspace{3ex} \forall (z', z'')\in \overset{p}{\underset{j=1}\bigcup} U_j,$ \\

\noindent and $\varepsilon(\delta)\downarrow 0$ when $\delta \downarrow 0$.

\vspace{2ex}

\noindent Since $\varepsilon(|z'|)\leq \varepsilon(\rho)$ for $z'\in B'(0, \, \rho)$, the previous estimate entails \\

\noindent $\frac{\displaystyle \sum\limits_{|\alpha|=k}\bigg|\frac{\frac{\partial^{\alpha}F_{k-1}}{\partial z^{'\alpha}}(0, z'')}{\alpha !}\bigg|^2 \, e^{-2\varphi (0, z'') - 2A|z''|^2}}{\displaystyle |\Lambda^r(ds)(0, z'')|^{2\frac{r+k}{r}}} \leq $ \\

\noindent $\displaystyle \leq \mbox{Const} \, \frac{2(r+k)}{\rho^{2(r+k)}} \,  e^{2(\varepsilon(\rho) + A \rho^2)} \,  \sup\limits_{(z', \, z'')\in U_j} \frac{|s(z', z'')|^{2r}\, (-\log |s(z', z'')|)^2}{|\Lambda^r(ds)(0, z'')|^{2\frac{r+k}{r}}} \,$

$$  \,  \int\limits_{z'\in B'(0, \, \rho)}\frac{||F_{k-1}(z', z'')||^2}{|s(z', z'')|^{2r} \, (-\log |s(z', z'')|)^2} \, d\lambda (z'),$$

\noindent for all $z''\in B''(0, \, \rho).$ A topological property of $Y$ ensures that there exists a nonnegative integer $N$ such that the covering $(U_j)_j$ of $Y_c$ can be chosen in such a way that $\#\{ j \, ; \, U_j \ni y\} \leq N.$ An integration with respect to $z''$ in the previous inequality, a summation on $j$, and obvious upper bounds yield \\

\noindent $\displaystyle \int\limits_{Y_c}\frac{|\nabla^k(J^kF_{k-1})|^2}{|\Lambda^r(ds)|^{2\frac{r+k}{r}}}dV_{Y, \, \omega} \leq C_{r, \, k}\, N\, M(c) \frac{1}{\rho^{2(r+k)}} e^{2(\varepsilon(\rho) + A\rho^2)} \int\limits_{\Omega'} \frac{||F_{k-1}||^2}{|s|^{2r}\, (-\log |s|)^2} dV_{X, \, \omega},$ \\

\noindent if $M(c)= \sup\limits_{(z', \, z'')\in \Omega'} \displaystyle \frac{|s(z', z'')|^{2r}\, (-\log |s(z', z'')|)^2}{|\Lambda^r(ds)(0, z'')|^{2\frac{r+k}{r}}}$ and $C_{r, \, k}=\mathrm{Const} \, 2(r+k).$   \\

\noindent The radius $\rho$ of the local holomorphic coordinate charts on which the subvariety $Y$ can be redressed is explicitly given by the following elementary lemma which is a refinement of the local inversion theorem to express the ``size'' of the ball on which we have a local diffeomorphism. 

\begin{Lem}\label{Lem:fonctionsimplicites} Let $E$ and $F$ be Banach spaces, $U$ an open subset of $E$, and 

\noindent $f:U \rightarrow F$ a $C^1$ map such that its differential map $df_a : E \rightarrow F$ at a point $a\in U$ is a bicontinuos isomorphism. 

 Then the open neighbourhood $V$ of $a,$ given by the local inversion theorem, on which $f$ is a diffeomorphism onto its image, contains the ball $B(a, \rho)$, where $\displaystyle \rho = \frac{1}{6(||df_a^{-1}||)(\sup\limits_{\xi\in U}||d^2 f_{\xi}||)}.$

\end{Lem}

  We leave the elementary proof of this lemma to the reader. It can be easily obtained from the proof of the local inversion theorem. Since the subvariety $Y$ is defined by the section $s \in H^0(X, \, E),$ we infer the explicit form of the weight function $\rho$ featuring in the statements of theorems \ref{The:principal1} and \ref{The:1}. Indeed, if $\theta : E_{|U} \rightarrow U \times \C^r$ is a trivialization of $E_{|U}$, and $(e_1, \dots ,\, e_r)$ the corresponding local holomorphic frame of $E_{|U},$ the restriction of $s$ to $U$ can be uniquely written as

\vspace{1ex}

  \hspace{6ex}  $s = \overset{r}{\underset{j=1}\sum} \sigma_j \otimes e_j,$  \hspace{3ex}  $\sigma_j \in {\cal O}(U).$

\vspace{1ex}

\noindent If $D$ is the Chern connection of the Hermitian holomorphic vector bundle $E,$ the operator $D$ can be written as

\vspace{1ex}

\hspace{6ex} $Ds \simeq_{\, \theta} d \, \sigma + A \wedge \sigma,$  

\vspace{1ex}

\noindent where $A=(a_{j \, k})$ is the matrix of $1$-forms representing the connection $D$ in the trivialization $\theta.$ Since the coefficients $a_{jk}$ of $A$ are locally bounded (by constants depending implicitly on $E$), lemma \ref{Lem:fonctionsimplicites} and the expression of $d$ in terms of $D$ show that the radius of the coordinate ball on which $Y$ can be redressed in a neighbourhood of a point $y \in Y$ is bounded below by \\

$\hspace{6ex} \displaystyle C \, \rho(y) = C \, \frac{1}{||Ds_y^{-1}|| \sup\limits_{\xi}(||D^2s_{\xi}|| + ||Ds_{\xi}||)},$  \\

\noindent the constant $C >0$ depending only on $E.$ This completes the proof of theorem \ref{The:1}.

\subsection{A Rauch-type comparison theorem}\label{subsection:untheoremedecomparaison}

  Recall that theorem \ref{The:principal1} was set on a K\"ahler manifold $(X, \, \omega)$. In order to get final estimates independent of the radius of local holomorphic coordinate balls of $X$, we prefer working on the tangent space to $X$ at a point instead of $X$ itself. The exponential map locally identifies $X$ to its tangent space. In order to estimate the deviation of the pull-back of $\omega$ to the tangent space from the standard Euclidian metric on this very tangent space, we need to establish a Riemannian geometric result related to the Rauch comparison theorem (see, for instance, [BC64], page 250). The proof of this result will be a slight reshaping of the proof of Rauch's theorem and will use the Jacobi vector fields theory and an elementary Gronwall-type lemma.

\vspace{2ex}

   Let $(M, g)$ be a complete Riemannian manifold, $m\in M$ an arbitrary point, and $\mathrm{exp}_m:T_mM \rightarrow M$, the exponential map at the point $m$. Let $\mathrm{Id}:= \mathrm{Id}_{T_mM}$ and, for an arbitrary point $x\in T_mM,$ consider the tangent linear map (or the differential) $T_x\mathrm{exp}_m : T_mM \rightarrow T_{\mathrm{exp}_m(x)}M$ of $\mathrm{exp}_m$ at the point $x$. We can identify $T_mM$ and $T_{\mathrm{exp}_m(x)}M$ via the isometry defined by parallel transport along the geodesic sprung from $x$. Our goal is to estimate

\vspace{1ex}

 $ ||T_x\mathrm{exp}_m - \mathrm{Id}|| $

\vspace{1ex}

\noindent in terms of $||x||,$ when $x$ ranges over the tangent space $T_mM$. Let $u\in T_mM, \,  ||u||=1,$ and $\gamma_u$ the geodesic sprung from $u$. We thus have

\vspace{1ex}

 $\gamma_u(0)=m$ et $\gamma_u(t)=\mathrm{exp}_m(tu)$,

\vspace{1ex}

\noindent for all $t$ in the definition interval of $\gamma_u$. Recall that a vector field $Y$ along the geodesic $\gamma_u$ is said to be a {\it Jacobi field} if it satisfies the second order differential equation \\

 $Y'' + R(\gamma_u', Y)\gamma_u'=0$, \\

\noindent where $R$ is the curvature tensor of $(M, \, g)$ defined as $R(X, Y)Z=\nabla_Y\nabla_XZ - \nabla_X\nabla_YZ + \nabla_{[X, Y]}Z$. It is a well-known fact that the differential of the exponential map is given by a Jacobi field. More precisely, for any $u, v\in T_mM$, we have the relation

\vspace{2ex}

 $(T_{tu}\mathrm{exp}_m)(tv)=Y(t)$, 

\vspace{2ex}

\noindent  where $Y$ is the unique Jacobi field along $\gamma_u$ such that $Y(0)=0$ and $Y'(0)=v$.

Assume now the sectional curvature of $(M, g)$ to be bounded, namely that there exists a constant $k>0$ such that

\vspace{1ex}

  $-k \leq K(p, P) \leq k $,

\vspace{1ex}

\noindent for every point $p\in M$ and every plane $P\subset T_pM$, where $K(p, P)$ stands for the sectional curvature of the plane $P$. So that we may estimate $||T_x\mathrm{exp}_m-\mathrm{Id}||$ we need estimate 

\vspace{2ex}

 $||(T_{tu}\mathrm{exp}_m)(tv)-\mathrm{Id}(tv)||=||Y(t)-Y'(0)t||$, \\

\noindent when $t$ ranges over $\R$. We need therefore an estimation of $Y$ which is known to satisfy a second order linear differential equation. The following elementary lemma, of Gronwall-type, provides the necessary estimation.

\begin{Lem}\label{Lem:Gronwall} Let $v:[0, T] \rightarrow \R$ be a $C^2$ function, $v \geq 0$, such that $v(0)=0, v'(0)=A$  and

\vspace{1ex}

 $-kv \leq v'' \leq kv$,  \hspace{2ex} on $[0, T]$,

\vspace{1ex}

\noindent  where $k>0$ is a constant. Then,

\vspace{1ex}

$A\frac{1}{\sqrt{k}}\sin (\sqrt{k}t) \leq v(t) \leq A\frac{1}{\sqrt{k}}\sinh(\sqrt{k}t),$ for all $t\in [0, T]$.

\end{Lem}

\noindent {\it Proof.} Let us first prove the right-hand inequality. Let $u$ be the solution to the Cauchy problem  $u''=ku$ with initial conditions $u(0)=0$ and $u'(0)=1$. Then, $u(t)=\frac{1}{\sqrt{k}}\sinh(\sqrt{k}t)$. In particular, $u\geq 0,$ and $u(t)=0$ if and only if $t=0$. The hypothesis shows that \\

$\frac{v''}{v}\leq k=\frac{u''}{u} \Longleftrightarrow (v'u-vu')'\leq 0 \Rightarrow v'u-vu'\leq 0$, \\

\noindent on $[0, T]$. This implies \\

$(\frac{v}{u})' \leq 0 \Rightarrow \frac{v(t)}{u(t)}\leq \frac{v}{u}(0_{+})$, \\

\noindent for all $t\in [0, T]$. Therefore, 

\vspace{2ex}

$v(t)\leq \frac{v}{u}(0_{+}) \, \frac{1}{\sqrt{k}}\sinh(\sqrt{k}t),$ \\

\noindent  for all $t\in [0, T]$. On the other hand, we see that

\vspace{2ex}

 $\frac{v}{u}(0_{+})=\lim\limits_{t\rightarrow 0}\frac{v(t)}{u(t)}=\lim\limits_{t\rightarrow 0}\frac{v'(t)}{u'(t)}=\frac{v'(0)}{u'(0)}=A,$

\vspace{2ex}

\noindent which proves the right-hand inequality. Let us now prove the left-hand inequality.

 Let $u$ be the solution to the Cauchy problem $u''=-ku$, with initial conditions $u(0)=0$ and $u'(0)=1$. Then, $u(t)=\frac{1}{\sqrt{k}}\sin(\sqrt{k}t)$. In particular, $u\geq 0,$ and $u(t)=0$ if and only if $t=0$. By hypothesis, we see that \\

$\frac{v''}{v} \geq -k=\frac{u''}{u} \Longleftrightarrow (v'u-vu')' \geq 0 \Longrightarrow v'u-vu'\geq 0$, \\

\noindent on $[0, T]$. This implies 

\vspace{2ex}

 $(\frac{v}{u})'\geq 0 \Rightarrow \frac{v(t)}{u(t)}\geq \frac{v}{u}(0_{+})$, \\

\noindent for all $t\in [0, T]$. Consequently,

\vspace{2ex}

 $v(t) \geq \frac{v}{u}(0_{+}) \, \frac{1}{\sqrt{k}}\sin(\sqrt{k}t),$ \\ 

\noindent for all $t\in [0, T]$. As before, $\frac{v}{u}(0_{+})=\frac{v'(0)}{u'(0)}=A$, which proves the left-hand inequality.   \hfill  $\Box$ \\

 We shall apply now this lemma to the components $Y_j$ of the Jacobi field $Y=(Y_1, \dots , \, Y_{2n})$ which are real functions satisfying $Y_j(0)=0,$ $Y'_j(0)=v_j,$ and $-k Y_j \leq Y''_j \leq k Y_j,$ for all $j=1, \dots , \, 2n,$ where $2n$ is the real dimension of the manifold $M$ and $v=(v_1, \dots , \, v_{2n})$ are the components of $v\in T_mM\simeq \R^{2n}.$ We get \\

  $|Y_j(t)-Y'_j(0)t|^2 \leq \bigg | \frac{\sinh (\sqrt{k}t)}{\sqrt{k}} - t \bigg |^2 \, |v_j|^2,$  \hspace{3ex} for $j=1, \dots , \, 2n,$ \\

\noindent if we also use that inequality $\sin x \leq x \leq \sinh x,$ for $x \geq 0.$ A summation on $j=1, \dots , \, 2n$ gives \\

$||Y(t)-Y'(0)t||\leq \bigg |\frac{\sinh(\sqrt{k}t)}{\sqrt{k}}-t\bigg |\, ||v||$, 

\vspace{1ex}

\noindent for all $t, v, u.$  We subsequently get, after dividing out by $t$, that 

\vspace{2ex}

 $||(T_{tu}\mathrm{exp}_m)(v)-\mathrm{Id}(v)||\leq \bigg |\frac{\sinh(\sqrt{k}t)}{\sqrt{k}t}-1 \bigg|\, ||v||$, \\

 $||T_{tu}\mathrm{exp}_m-\mathrm{Id}|| \leq \bigg |\frac{\sinh(\sqrt{k}t)}{\sqrt{k}t}-1\bigg |$, \\

\noindent for all $t, u$. If we set $x=tu,$ we find

\vspace{2ex}

$||T_x\mathrm{exp}_m-\mathrm{Id}||\leq \bigg |\frac{\sinh(\sqrt{k}||x||)}{\sqrt{k}||x||}-1 \bigg |$, for all $x\in T_mM$.

\vspace{2ex}

\noindent  Since $\sinh x\geq x,$ for all $x\geq 0$, the absolute value is superfluous in the right-hand term. We have thus proved the following.

\begin{Prop}\label{Prop:Rauch}

 If there exists a constant $k>0$ such that

\vspace{1ex}

 $ -k \leq K(p, P) \leq k$,

\vspace{1ex}

\noindent for every point $p\in M$ and every plane $P\subset T_pM$, then

\vspace{1ex}

$\displaystyle ||T_x\mathrm{exp}_m-\mathrm{Id}||\leq \frac{\sinh(\sqrt{k}||x||)}{\sqrt{k}||x||}-1,$ \hspace{2ex} for all $x\in T_mM$.

\end{Prop}

\noindent {\bf Remark.} The Rauch comparison theorem estimates $||T_x\mathrm{exp}_m||$. The above proposition estimates the distance between $T_x\mathrm{exp}_m$ and $T_0\mathrm{exp}_m=\mathrm{Id}$. The latter is therefore slightly more general.

\subsection{Final estimate}\label{subsection:estimationfinale}

In order to complete the proof of theorem \ref{The:principal1}, it remains to get a uniform control of $\displaystyle\int\limits_{Y_c}\frac{|\nabla^k(J^kF_{k-1})|^2}{|\Lambda^r(ds)|^{2\frac{r+k}{r}}} dV_{Y, \, \omega}$  \hspace{1ex} (see the end of section \ref{subsection:demonstrationdutheoremeprincipal}).

  Fix a point $y_0\in Y \subset X$, and let $\Phi := \exp_{y_0} : T_{y_0}X \rightarrow X$ be the exponential map. The K\"ahler metric $\omega$ on the weakly pseudoconvex manifold $X$ can be made complete by a standard well-known procedure. We may therefore assume, without loss of generality, that the exponential map is defined on the whole tangent space. Let $\omega_0$ be the standard K\"ahler metric on the Euclidian space $T_{y_0}X \simeq \C^n$. Our first goal in this section is to find an explicit formula for the radius of the ball in the tangent space $T_{y_0}X$ on which the two metrics $\Phi^{\star}\omega$ and $\omega_0$ can be compared. Let us set \\

$(0. 6. 1)$  \, $r(y_0):= \sup\{r>0 \, ; \sup\limits_{\stackrel{x\in B(y_0, \, r)} {0\leq l \leq m}}r^{2+l}\, ||\nabla^l\Theta(T_X)(x)||<10^{-2a}\}$,\\

\noindent where $a>0$ is a constant to be specified later, and $\nabla^l\Theta(T_X)$ stands for the $l^{\mbox{th}}$ order derivative of the curvature tensor $\Theta(T_X)$ viewed as a section of the $C^{\infty}$ bundle $\Lambda^{1, \, 1}T^{\star}_X \otimes \mathrm{Hom}(T_X, \, T_X).$ Locally, this boils down to deriving the coefficients of $\Theta(T_X).$ In particular, we get

\vspace{1ex}

$\sup\limits_{x\in B(y_0, \, r(y_0))}||\Theta(T_X)|| \leq \frac{10^{-2a}}{r(y_0)^2}:=k$, 

\vspace{1ex}

\noindent and hence the following bounds for the sectional curvature of the manifold $X$ :

\vspace{1ex}

\hspace{6ex} $ -k \leq K(p, P) \leq k,$ 

\vspace{1ex}

\noindent for all $p\in B(y_0, r(y_0))$, and all planes $P \subset T_{y_0}X$ in the tangent space at $y_0$ to $X.$

 This shows that the hypothesis of proposition \ref{Prop:Rauch} is fulfilled in the ball $B(y_0, r(y_0))$. Then we get \\

 $(\star)$  \hspace{2ex}  $\displaystyle ||T_v\mathrm{exp}_{y_0}-\mathrm{Id}|| \leq \frac{\sinh(\sqrt{k}||v||)}{\sqrt{k}||v||}-1,$ \\

\noindent  for all $v\in T_{y_0}X,$ such that $||v||<r(y_0)$. If $||T_v\mathrm{exp}_{y_0}-\mathrm{Id}||<1$, the map $T_v\mathrm{exp}_{y_0}$ is invertible. Consequently, $\mathrm{exp}_{y_0}$ is an immersion on $B(0, r(y_0)) \subset T_{y_0}X,$ if  $\frac{\sinh(\sqrt{k}||v||)}{\sqrt{k}||v||}<2$ for all $v$ such that $||v||<r(y_0).$ To achieve this, it is enough to have \\

 $(1)$ \,  $\frac{\sinh(10^{-a})}{10^{-a}}<2. $ \\

 On the other hand, we need a value of the constant $a$ such that we may have the bounds \\

 $(\star\star)$ \hspace{2ex} $\frac{1}{2}\omega_0 \leq \mathrm{exp}_{y_0}^{\star}\omega \leq 2\omega_0$, \hspace{2ex} on the ball $B(0, r(y_0))$ in $T_{y_0}X$.

\vspace{2ex}

\noindent In order to have these bounds, it is enough to have

\vspace{2ex}

 $\frac{1}{2} \leq ||T_v\mathrm{exp}_{y_0}||\leq 2$,

\vspace{2ex}

\noindent for all $v\in T_{y_0}X$ such that $||v||<r(y_0)$. We thus infer from $(\star)$ that

\vspace{2ex}

$2-\frac{\sinh(\sqrt{k}||v||)}{\sqrt{k}||v||} \leq ||T_v\mathrm{exp}_{y_0}|| \leq \frac{\sinh(\sqrt{k}||v||)}{\sqrt{k}||v||}$,

\vspace{2ex}

\noindent for all $v\in T_{y_0}X, ||v||<r(y_0)$. This shows that it is enough to have $\frac{\sinh(\sqrt{k}||v||)}{\sqrt{k}||v||} \leq \frac{3}{2}$, for all $v$ such that $||v||<r(y_0)=\frac{10^{-a}}{\sqrt{k}}$. The bounds $(\star \star)$ are therefore guaranteed as soon as the constant $a$ satisfies the inequality \\

 $(2)$ \,  $\frac{\sinh(10^{-a})}{10^{-a}} \leq \frac{3}{2}.$ 

\vspace{2ex}

\noindent In short, we have proved the following. 

\begin{Lem}\label{Lem:exponentielle} For a choice of the constant $a> 0$ satisfying inequality $(2)$, and for $r(y_0)$ defined by relation $(0. 6. 1),$ the exponential map $\Phi=\exp_{y_0}$ is an immersion and the bounds $(\star\star)$ hold on the ball $B(0, r(y_0))$ in the tangent space $T_{y_0}X.$

\end{Lem}

 Lemma \ref{Lem:fonctionsimplicites} shows that there exist local holomorphic coordinates $\zeta=(\zeta', \, \zeta''), \, \zeta'=(\zeta_1, \dots , \, \zeta_r),$  $\zeta''=(\zeta_{r+1}, \dots , \, \zeta_n)$ on the ball $B(0, \, r) \subset T_{y_0}X$ such that the subvariety $\Phi^{-1}(Y\cap B(y_0, \, r))\subset B(0, \, r)$ is defined by the equations $\zeta'=0$, for the following radius \\

$\displaystyle r= \rho(y_0) = \frac{1}{6 \, ||Ds_{y_0}^{-1}||_{\omega_0}\, \sup\limits_{\xi}(||(D^2s_{\xi}||_{\omega_0}+ ||Ds_{\xi}||_{\omega_0})}.$

\noindent Moreover, the bounds $(\star\star)$ imply \\

$\displaystyle r \geq \frac{1}{24\, ||Ds_{y_0}^{-1}||_{\omega}\, \sup\limits_{\xi}(||D^2s_{\xi}||_{\omega} + ||Ds_{\xi}||_{\omega})}:=r_0 (y_0). $ \\

\noindent In the above expressions all $\sup\limits_{\xi}$ are computed for $\xi\in B(y_0, r(y_0))$. Let us set from now on: 

\vspace{1ex}

 $(0. 6. 2)$ \hspace{6ex}  $r_1(y_0) = \min(r(y_0), \, r_0(y_0))$.

\vspace{1ex}

 Recall that $F_{k-1}\in H^0(X, \, \Lambda^nT^{\star}_X \otimes L)$ is the extension of the jet $f\in H^0(X, \, \Lambda^nT^{\star}_X \otimes L \otimes {\cal O}_X/{\cal I}_Y^{k+1})$ to order $k-1$ given by the induction hypothesis of theorem \ref{The:principal1} (see the beginning of \ref{subsection:demonstrationdutheoremeprincipal}). The holomorphic line bundle $L':= \Lambda^nT^{\star}_X \otimes L$ is equipped with a $C^{\infty}$ Hermitian metric $h$. Let us consider the $C^{\infty}$ line bundle $\Phi^{\star}L'$ equipped with the metric $\phi^{\star}h$ and the section $\Phi^{\star} F_{k-1} \in C^{\infty}(T_{y_0}X, \, \Phi^{\star}L').$

 Let $J_X \in \mathrm{End}(T_X)$ be the complex structure of the manifold $X$ and $J:= \Phi^{\star}J_X$ the almost complex structure induced on $T_{y_0}X.$ If $J_0$ is the canonical complex structure of $T_{y_0} \simeq \C^n,$ the map $\Phi$ is not $(J_0,  J_X)$-holomorphic, but it certainly is $(J,  J_X)$-holomorphic. If $i\Theta(L')$ is the curvature form (of type $(1, \, 1)$) of $(L', \, h),$ $\Phi^{\star}(i\Theta(L'))$ is a type $(1, \, 1)$-form for $J$ on $T_{y_0}X.$

\begin{Lem}\label{Lem:potentiel1} There exists a real function $\tilde{\varphi} \in C^{\infty}$ on the ball $B= B(0, \, r_1(y_0))$ in the tangent space $T_{y_0}X$ such that $i\partial_J \bar{\partial}_J \tilde{\varphi}=\Phi^{\star}(i\Theta(L'))$ and

\vspace{1ex}

$\sup\limits_B |\tilde{\varphi}| \leq C \sup\limits_B ||\Phi^{\star}(i\Theta(L'))||,$ 

\vspace{1ex}

\noindent where $C>0$ is a constant depending only on $r_1(y_0).$

\end{Lem}

\noindent {\it Proof.} With respect to real coordinates $x_1, \dots , \, x_{2n}$ on $B,$ the real $d$-closed $2$-form $\Phi^{\star}(i\Theta(L'))$ can be written as $\Phi^{\star}(i\Theta(L')) = \sum\limits_{i < j} v_{ij} \, dx_i \wedge dx_j,$ with functions $v_{ij} \in C^{\infty}(B).$ The Poincar\'e lemma gives the explicit formula: \\

$U(x)= \sum\limits_{i < j} \bigg(\displaystyle \int_0^1 t \, v_{ij}(tx) \, dt \bigg ) (x_i \, dx_j - x_j \, dx_i),$ \\

\noindent for a $C^{\infty}$ solution of the equation $dU=\Phi^{\star}(i\Theta(L'))$ on $B$. We see then that

\vspace{1ex}

$||U||_{L^{\infty}(B)} \leq C_1 \, ||\Phi^{\star}(i\Theta(L'))||_{L^{\infty}(B)},$

\vspace{1ex}

\noindent with a constant $C_1 > 0$ depending only on the radius of $B$. With respect to the almost complex structure $J,$ the real $1$-form  $U$ decomposes as $U = U^{1, 0} + U^{0, 1},$ with $U^{0, 1}= \overline{ U^{1, 0}}.$ Then $dU= \partial_J U^{0, 1} + \overline{\partial_J U^{0, 1}},$ since $dU$ is of type $(1, \, 1)$ for $J$. The almost complex structure $J$ is integrable as the inverse image of an integrable almost complex structure. Let $(z_1, \dots , \, z_n)$ be $J$-holomorphic complex coordinates centred at $0$ on a neighbourhood of the ball $B \subset T_{y_0}X.$ We thus have $\bar{\partial}_JU^{0, \, 1} = 0$ on $B$. The bounds $(\star\star),$ relating the metrics $\omega$ and $\omega_0,$ allow us to assume that the ball $B$ is $J$-pseudoconvex (if not so, we multiply the radius $r_1(y_0)$ by a fixed constant). Since for an integrable almost complex structure we have the same formalism as for an complex analytic structure, a classical result on the solvability of the $\bar{\partial}$ operator on bounded strictly pseudoconvex domains with a $C^2$ boundary in $\C^n$ (see, for instance, [HL84], theorem 2.3.5.), yields the existence of a constant $C_2>0$ depending only on the radius of the ball $B$, and of a solution to the equation $\bar{\partial}_J v = U^{0, 1}$ on $B$ obtained by an explicit integral formula, such that

\vspace{1ex}

$||v||_{L^{\infty}(B)}  \leq C_2 \, ||U^{0, 1}||_{L^{\infty}(B)} \leq 2  C_2 \, ||U||_{L^{\infty}(B)}.$

\vspace{1ex}

\noindent Then $\tilde{\varphi}:= i(\bar{v} - v)$ is the function we were looking for.    \hfill $\Box$   \\

  Since $\phi$ is an immersion on $B(0, r_1(y_0))$, there exists a neighbourhood $V \subset B(0, \, r_1(y_0))$ of $0$ such that  $\phi$ is a diffeomorphism of $V$ onto a neighbourhood $U$ of $y_0$ in $X$. Let $\psi : U \rightarrow V$ be the inverse diffeomorphism. In a local trivialization of $L'$ in a neighbourhood of $y_0,$ the section $F_{k-1}$ can be written as $F_{k-1} = u \otimes e,$ for a local holomorphic frame $e$. The function $v=u \circ \Phi$ is then $C^{\infty}$ on $V$, and $u$ being holomorphic implies: $\bar{\partial}(v\circ \psi)=0.$ If $z=(z_1, \dots , \, z_n)$ is a system of local holomorphic coordinates on $U$, this means that $v$ is a solution to the following elliptic system

\vspace{1ex}

 $(\star\star\star)$ \hspace{2ex} $\displaystyle \sum\limits_j \frac{\partial v}{\partial \zeta_j}\circ \psi \, \frac{\partial \psi_j}{\partial\bar{z}_k} + \sum\limits_j \frac{\partial v}{\partial\bar{\zeta}_j}\circ \psi \, \frac{\partial \bar{\psi}_j}{\partial \bar{z}_k} = 0,$  \, $k=1, \dots , \, n.$

\vspace{1ex}

\noindent  Let us remind now a standard differential operator theory result. G\aa rding's lemma controls the growth of the derivatives of a solution to an elliptic equation in terms of the growth of this very solution. This lemma plays the role of Cauchy's inequalities in the nonholomorphic case. Let $H_j^{loc}$ be the Sobolev space of locally $L^2$ functions whose all derivatives in the sense of distributions up to order $j$ are still locally $L^2,$ and let $||\,\,  ||_j$ be its Sobolev norm. We refer for the details to [Agm65] (lemma 6.1 and theorems 6.2-6.7, pages 53-67).

\begin{The}(theorem $6.5$ in [Agm65]) Let $\Omega$ be an open subset of $\R^n$, and $A_1(x, D),$ $\dots , A_N(x, D)$ differential operators of respective orders $m_1, \dots , \, m_N,$  with coefficients $a_{\alpha}^i \in C^{\infty},$ which make up an elliptic system in $\Omega$. Let $u\in L^2_{\mathrm{loc}}(\Omega)$ such that $A_i^{\star}u\in H_{k_i}^{\mathrm{loc}}(\Omega)$, for all $i=1, \dots , N$. 

  If $j:=\min (m_1+k_1, \dots , m_N+k_N)$, then $u\in H_j^{loc}(\Omega).$ In addition, for all $\Omega' \subset\subset  \Omega$, there exists $\gamma = \gamma(A_i, \Omega', \Omega)$ such that

\vspace{1ex}

 $||u||_{j, \, \Omega'} \leq \gamma \, (\overset{N}{\underset{i=1}{\sum}}|A_i^{\star}u|_{k_i, \, \Omega} + ||u||_{0, \, \Omega})$,

\vspace{1ex}

\noindent where $\gamma = \mathrm{Const} \cdot p\cdot N \cdot K \cdot M$, $\mathrm{Const}$ is a universal constant, $p=p(n, l)=\mathrm{card}\{\alpha\in \mathbb{N}^n | |\alpha|=l\}$, $K=\sup\limits_{\xi\in \Omega', \, |\alpha|\leq l, i}|da_{\alpha}^i(\xi)|$, $M=\sup\limits_{x\in \Omega', \, |\alpha|\leq l, i}|a_{\alpha}^i(x)|$. 

\end{The}

 The actual dependence of the constant $\gamma$ on the data is not explicit in [Agm65], but it can be easily inferred from the proofs given there to theorems $6.2-6.7$. Also, there is a slightly more general statement there in which the coefficients of the operators $A_i(x, D)$ are only assumed to be ``$s$-smooth''.  

  Since $v$ is a solution to the elliptic system $(\star\star\star)$, the previous theorem shows that we have the estimate \\

 $\displaystyle \sup\limits_{||\zeta''|| \leq  \frac{1}{2} r_1(y_0)} \sum\limits_{|\alpha|=k} \bigg | \frac{\partial^{\alpha}v}{\partial \zeta^{'\alpha}}(0, \, \zeta'')\bigg |^2  \leq \gamma_k \,  \int_{B(0, \, r_1(y_0))} |v(\zeta', \, \zeta'')|^2 \, d\, \lambda(\zeta', \, \zeta''),$ \\

\noindent where $\gamma_k=\mathrm{Const} \cdot p_k\cdot \max(\sup\limits_{\xi\in U}||d_{\xi}\psi||, \, \sup\limits_{\xi\in U}||d_{\xi}^2\psi||)$, $p_k = \mathrm{Card}\{\alpha | \, |\alpha|=k\}$ and $\mathrm{Const}$ is a universal constant. For the following norms computed in the Hermitian vector bundle $(\Phi^{\star}L', \, \Phi^{\star} h)$, equipped with the local weight $\tilde{\varphi}$,  \\

$\bigg |\bigg | \frac{\partial^{\alpha}v}{\partial \zeta^{'\alpha}}(0, \, \zeta'')\bigg |\bigg |^2 = \bigg | \frac{\partial^{\alpha}v}{\partial \zeta^{'\alpha}}(0, \, \zeta'')\bigg |^2 \, e^{-2\tilde{\varphi}(0, \, \zeta'')}, $ \hspace{2ex} $||v(\zeta', \, \zeta'')||^2 = |v(\zeta', \, \zeta'')|^2 \, e^{-2\tilde{\varphi}(\zeta', \, \zeta'')}, $ \\

\noindent we get the estimate \\

$\displaystyle \int\limits_{||\zeta''|| \leq  \frac{1}{2} r_1(y_0)} \sum\limits_{|\alpha|=k} \bigg |\bigg | \frac{\partial^{\alpha}v}{\partial \zeta^{'\alpha}}(0, \, \zeta'')\bigg |\bigg |^2 d\, \zeta''  \leq $  \\

\hspace{3cm}  $\leq \displaystyle \gamma_k \,  \int_{B(0, \, r_1(y_0))} ||v(\zeta', \, \zeta'')||^2 \, e^{2(\tilde{\varphi}(\zeta', \zeta'') - \tilde{\varphi}(0, \zeta''))} \, d\, \lambda(\zeta', \, \zeta''),$ \\

\noindent and also, thanks to lemma \ref{Lem:potentiel1},   \\

\noindent $(3)$ \, $\displaystyle \int\limits_{||\zeta''|| \leq \frac{1}{2} r_1(y_0)} \sum\limits_{|\alpha|=k} \bigg |\bigg | \frac{\partial^{\alpha}v}{\partial \zeta^{'\alpha}}(0, \, \zeta'')\bigg |\bigg |^2 d\, \zeta'' \leq \gamma_k \, C_{L'} \,  \int_{B(0, \, r_1(y_0))} ||v(\zeta', \, \zeta'')||^2 \, d\, \lambda(\zeta', \, \zeta''), $ \\

\noindent where the constant $\displaystyle C_{L'}:= e^{2C \sup\limits_U ||i\Theta(L')||}$ depends only on the growth of the curvature of $L'.$

  It remains to infer from the estimate $(3)$ for $v$ an analogous estimate for $u$. If $z$ is the variable on $U \subset X$, and $\zeta$ is the variable on $V \subset T_{y_0}X$, the change of variable $\zeta = \psi(z)$ implies the following estimate for $u$   \\

$(4)$  \hspace{2ex}  $||u||^2_{k, \, U' \cap Y} \leq \tilde{\gamma}_k \, C_{L'} \,  ||u||^2_{0, \, U}$,  \, \, $U' \subset \subset U,$

\vspace{2ex}

\noindent where $\tilde{\gamma}_k = \mathrm{Const}\cdot p_k \cdot \sup\limits_{\stackrel{1\leq l \leq k}{\xi\in U}}||d_{\xi}^l\psi||$, $\mathrm{Const}$ being a universal constant.

   Proposition 3 already gave an estimate for the norm of the differential map of $\phi$, and implicitly for the differential map of $\psi$. The formula for $\tilde{\gamma}_k$ would also require an estimation of the growth of the differentials of order $\leq k$ of $\psi$. It is clear that $\sup\limits_{\stackrel{1\leq l \leq k}{\xi\in U}}||d_{\xi}^l\psi||$ is bounded above by a constant depending only on the radius $r_1(y_0)$ of the ball on which we are working. These are standard computations that can well be left to the reader. 

 We are now in a position to conclude that the constant $C_r^{(k)}$ in the statement of theorem \ref{The:principal1} depends only on $r$, on $k,$ on $E$, and on $\sup\limits_{\Omega}||i\Theta(L)||.$

\subsection{The case of a singular subvariety}\label{subsection:demonstrationdutheoremevarprincipal}

 A standard argument shows that the restriction imposed at the beginning of section \ref{subsection:demonstrationdutheoremeprincipal} on the singular set $\Sigma = \{ s =0, \, \Lambda^r(ds) = 0 \}$ of $Y$ to be empty is superfluous. Indeed, since the section $s \in H^0(X, \, E)$ is assumed to be generically transverse to the zero section, we can find a complex hypersurface $Z \subset X$ such that $\Sigma \subset \overline{Y}\cap Z \subsetneq \overline{Y}.$ If the ambient manifold $X$ is Stein, the complementary of $Z$ is still Stein. We can therefore apply theorem \ref{The:principal1} to the Stein manifold $X\setminus Z$ and use lemma \ref{Lem:dbarextension} to extend the $L^2$ estimates across $Z.$ In the general case of a weakly pseudoconvex ambient manifold $X$, we can apply the Stein case on coordinate balls $U_j$ to construct local holomorphic extensions $\tilde{f}_j$ of the jet $f$ satisfying estimates $\int_{U_j}|\tilde{f}_j|^2 |s|^{-2r}(-\log|s|)^{-2} \, dV < +\infty,$ and we set $\tilde{f}_{\infty} = \sum\limits_j \theta_j \, \tilde{f}_j,$ for a partition of unity $(\theta_j)_j.$ \hfill $\Box$

\vspace{3ex}

\noindent {\Large\bf References}

\noindent [Agm65] \, S. Agmon --- {\it Lectures on Elliptic Boundary Value Problems} --- Van Nostrand, Prineton, 1965.

\vspace{1ex}

\noindent [BC64] \, R. Bishop, R. J. Crittenden --- {\it Geometry of Manifolds} --- Academic Press, 1964.

\vspace{1ex}

\noindent [Dem 82] \, J.-P. Demailly --- {\it Estimations $L^2$ pour l'op\'erateur $\bar{\partial}$ d'un fibr\'e vectoriel holomorphe semi-positif au-dessus d'une vari\'et\'e k\"ahl\'erienne compl\`ete} --- Ann. Sci. Ecole Norm. Sup. {\bf 15}  (1982) 457-511.

\vspace{1ex}

\noindent [Dem 97] \, J.-P. Demailly --- {\it Complex Analytic and Algebraic Geometry}---http://www-fourier.ujf-grenoble.fr/~demailly/manuscripts/agbook.ps.gz

\vspace{1ex}

\noindent [Dem99] \, J.-P. Demailly --- {\it On the Ohsawa-Takegoshi-Manivel $L^2$ Extension Theorem} --- Article en l'honneur de Pierre Lelong \`a l'occasion de son 85\`eme anniversaire.

\vspace{1ex}

\vspace{1ex}

\noindent [HL84] \, G. Henkin, J. Leiterer --- {\it Theory of Functions on Complex Manifolds} --- Birkh\"auser Verlag, Basel, Boston, Stuttgart, 1984.

\vspace{1ex}

\noindent [H\"or65] \, L. H\"ormander --- {\it $L^2$ Estimates and Existence Theorems for the $\bar{\partial}$ Operator} --- Acta Math. {\bf 113} (1965) 89-152.

\vspace{1ex}

\noindent [H\"or66] \, L. H\"ormander --- {\it An Introduction to Complex Analysis in Several Variables} --- 1st edition, Elsevier Science Pub., New York, 1966, 3rd revised edition, North-Holland math. library, Vol 7, Amsterdam (1990).

\vspace{1ex}

\noindent [Man93] \, L. Manivel --- {\it Un th\'eor\`eme de prolongement $L^2$ de sections holomorphes d'un fibr\'e hermitien}--- Math. Zeitschrift {\bf 212} (1993) 107-122.

\vspace{1ex}

\noindent [OT87] \, T. Ohsawa, K. Takegoshi --- {\it On The Extension of $L^2$ Holomorphic Functions}--- Math. Zeitschrift {\bf 195} (1987) 197-204.

\vspace{1ex}

\noindent [Ohs88] \, T. Ohsawa --- {\it On the Extension of $L^2$ Holomorphic Functions, II}---Publ. RIMS, Kyoto Univ. {\bf 24} (1988) 265-275.

\vspace{1ex}

\noindent [Ohs94] \, T. Ohsawa --- {\it On the Extension of $L^2$ Holomorphic Functions, IV : A New Density Conept}--- Mabuchi, T (ed.) et al., Geometry and Analysis on Complex Manifolds. Festschrift for Professor S. Kobayashi's 60th birthday. Singapore:World Scientific, (1994) 157-170.

\vspace{1ex}

\noindent [Ohs95] \, T. Ohsawa --- {\it On the Extension of $L^2$ Holomorphic Functions, III: Negligible Weights} --- Math. Zeitschrift {\bf 219} (1995) 215-225.

\vspace{1ex}

\vspace{3ex}

\noindent Dan Popovici

\noindent Mathematics Institute

\noindent University of Warwick

\noindent Coventry CV4 7AL

\noindent United Kingdom

\noindent E-mail: popovici@maths.warwick.ac.uk

\end{document}